\documentclass[12pt]{amsart}
\usepackage{epsf,amssymb}
\usepackage{graphicx,epsf}

\theoremstyle{plain}
\newtheorem{thm}{Theorem}[section]
\newtheorem{lemma}[thm]{Lemma}
\newtheorem{prop}[thm]{Proposition}
\newtheorem{cor}[thm]{Corollary}
\theoremstyle{definition}
\newtheorem*{rmk}{Remark}
\newtheorem*{example}{Example}
\newtheorem{defn}[thm]{Definition}

\newcommand{\C}{{\mathbb{C}}}

\newcommand{\R}{{\mathbb{R}}}

\newcommand{\Z}{{\mathbb{Z}}}
\newcommand{\Q}{{\mathbb{Q}}}

\newcommand{\sig}{\sigma}
\newcommand{\bdy}{{\partial}}
\newcommand{\Sig}{{\Sigma}}
\newcommand{\ori}{\mathop{\mathit or}\nolimits}
\newcommand{\st}{\mathop{\rm st}\nolimits}
\newcommand{\Sym}{\mathop{\rm Sym}\nolimits}
\newcommand{\Min}{\mathop{\rm Min}}
\newcommand{\cal}{\mathcal}

\newcommand{\cE}{{\cal E}}
\newcommand{\cF}{{\cal F}}
\newcommand{\cG}{{\cal G}}

\newcommand{\cL}{{\cal L}}
\newcommand{\cP}{{\cal P}}
\newcommand{\cA}{{\cal A}}

\newcommand{\cM}{{\cal M}}

\newcommand{\ol}{\overline}

\newcommand{\la}{\langle}
\newcommand{\ra}{\rangle}
\newcommand{\lra}{\longrightarrow}
\newcommand{\wt}{\widetilde}
\newcommand{\Hom}{\mathrm{Hom}}

\newcommand{\m}{\mathfrak{m}}
\newcommand{\coker}{\mathop{\mathrm{coker}}}

\newcommand{\Span}{\mathop{\mathrm{span}}}
\newcommand{\Hilb}{\mathop{\rm Hilb}}
\newcommand{\Ru}{\underline{\R}}
\newcommand{\udot}{{\scriptscriptstyle \bullet}}
\newcommand{\oudot}{{\raisebox{2pt}[0pt]{$\udot$}}} 

\def\im{\mathop{\rm Im}\nolimits}
\newcommand{\Tor}{\mathop{\mathrm{Tor}}\nolimits}
\begin{document}

\title{Remarks on the combinatorial intersection cohomology of fans}

\date{\today}
\author{Tom Braden}
\address{University of Massachusetts, Amherst}
\email{braden@math.umass.edu}

\begin{abstract} 
We review the theory of combinatorial
intersection cohomology of fans developed by
Barthel-Brasselet-Fieseler-Kaup, Bressler-Lunts, and Karu.  This
theory gives a substitute for the intersection cohomology of toric
varieties which has all the expected formal properties but makes sense
even for non-rational fans, which do not define a toric variety.  As a
result, a number of interesting results on the toric $g$ and $h$
polynomials have been extended from rational polytopes to general
polytopes.  We present explicit complexes computing the combinatorial
IH in degrees one and two; the degree two complex gives the rigidity
complex previously used by Kalai to study $g_2$.  We present several
new results which follow from these methods, as well as previously
unpublished proofs of Kalai that $g_k(P) = 0$ implies $g_k(P^*) = 0$
and $g_{k+1}(P) = 0$.

\end{abstract}

\dedicatory{
For Bob MacPherson on his sixtieth birthday
}

\thanks{The author was supported in part by NSF grant DMS-0201823}

\maketitle

\newcommand{\gt}{\tilde{g}}
\newcommand{\zcone}{{o}}

\addtolength{\parskip}{4pt plus 2pt minus 2pt}

For a $d$-dimensional convex polytope $P$, Stanley \cite{St87} defined
a polynomial invariant $h(P, t) = \sum^d_{k = 0} h_k(P)t^k$ of its
face lattice, which is usually called the ``generalized"
or ``toric" $h$-polynomial of $P$.  It is ``generalized" in that
it extends a previous definition 
from simplicial polytopes to general polytopes, while
the adjective ``toric" refers to the fact that if $P$ is a rational
polytope (meaning that all its vertices have all coordinates in $\Q$), then
the coefficients of $h(P, t)$ are intersection cohomology 
Betti numbers of an associated projective toric variety $X_P$:
\begin{equation}\label{h = IH} h_k(P) = \dim_\R IH^{2k}(X_P; \R).
\end{equation}
In fact, these Betti numbers had been computed independently by several
people, including Robert MacPherson, and these calculations
inspired Stanley's definition (although proofs \cite{DL,Fie} of 
\eqref{h = IH} did not appear in print until several years after Stanley's 
definition of the $h$-polynomial).

This connection between the topology of toric varieties and combinatorics
of polytopes was then used to prove a number of interesting
relations among these invariants \cite{Bay98,BM,St92}; these results were
therefore known only for rational polytopes.  It was 
believed, however, that they should hold for general polytopes, and in 
a few cases (simplicial polytopes \cite{McM} and low degree terms
\cite{Kal:rigidity}) non-toric proofs were found.
 
The recent papers \cite{BBFK,BrLu,Ka} have settled this question, 
by defining groups which can substitute for intersection cohomology
in \eqref{h = IH}, but which are defined whether or not $P$ is rational.
In fact, they accomplish more: they define a theory
of sheaves purely in terms of the linear structure of $P$ which
completely captures the structure of torus-equivariant constructible sheaves 
on the toric
variety $X_P$ when $P$ is rational.  The maps and other relations between intersection cohomology
groups which were used to prove identities and inequalities
among the $h$-numbers \cite{BM,St92} have direct analogues in the
new theory.  Moreover, the deep theorems
from algebraic geometry that were used to prove \eqref{h = IH} --- 
the Decomposition Theorem and the Hard Lefschetz
Theorem --- now have proofs in this combinatorial setting.  The result is a 
powerful, self-contained theory which can answer many questions about the
combinatorics of convex polytopes.  

This paper is meant as a guided introduction to the theory of combinatorial 
intersection cohomology and some of its 
applications.  In Section 1 we describe the toric $g$- and $h$-numbers
and the main results about them.  In Section 2 we 
review the construction of the combinatorial intersection cohomology sheaves 
and their main properties.  We describe what their graded pieces look like
in degrees up to two; in degree two we recover a
chain complex related to rigidity of frameworks which Kalai used as a 
substitute for the intersection cohomology group $IH^4(X_P;\R)$.

In Section 3 we present two applications of a result of \cite{BBFK} on
``quasi-convex'' fans which have not appeared before.  
First, we answer a question of Stanley \cite{St87} 
regarding whether local contributions to $h(P,t)$ at a facet have nonnegative 
coefficients.  Second, we present an inequality which generalizes Kalai's monotonicity 
with $t$ specialized to $1$.  
In Section 4 we describe how Stanley's convolution 
identity which relates the $g$-numbers of $P$ and the polar polytope $P^*$
can be ``lifted'' to exact sequences, and relate this to the Koszul duality
constructed in \cite{BraLu}. Finally, in an appendix we present Kalai's 
proofs of two interesting applications of monotonicity.   

Although our discussion is developed entirely in the setting of convex 
geometry, in many places we point out the connections with the topology of
toric varieties.  Although this is no longer necessary to understand the 
proofs, we hope that the reader will appreciate seeing these connections 
as a way to motivate the constructions and suggest new applications.  

Thanks are due to Gottfried Barthel, Jean-Paul Brasselet, Karl-Heinz Fieseler,
Ludger Kaup, Valery Lunts, and Gil Kalai for many interesting conversations
over the years, to Chris McDaniel for helpful comments on a draft of this
article, and to Bob MacPherson for introducing me to the beautiful interplay
between geometry and combinatorics. 

\section{The $g$- and $h$-polynomials} 
\subsection{Simplicial polytopes}
We start with a discussion of the situation for simplicial polytopes.
When $P$ is simplicial, the $h$-polynomial $h(P, t) = \sum^d_{k = 0} h_k(P)t^k$
is defined by
\begin{equation}\label{simplicial h formula}
h(P,t) = (t - 1)^d + f_0(t-1)^{d-1} + \dots + f_{d-1},
\end{equation}
where the face number $f_k = f_k(P)$ is the number of $k$-dimensional faces
of $P$.  The transformation taking $\{f_k\}$ to $\{h_k\}$ 
is invertible, so that the $h$-numbers determine the face numbers.  

The $h$-numbers are not
independent --- they satisfy 
the {Dehn-Sommerville relations}: 
\[h_k(P) = h_{d-k}(P)\; \text{for
all }\; 1 \le k \le d.\]  
As a result, if we define an auxiliary polynomial 
$g(P, t) = \sum^{\lfloor d/2\rfloor}_{k = 0} g_k(P)t^k$
by \begin{equation}\label{g from h}
g_k(P) = h_k(P) - h_{k-1}(P),\; 0\le k \le d/2,
\end{equation}
where we put $h_{-1}(P)=0$, then
$g(P, t)$ contains all the information 
of $h(P, t)$, and thus all the information about 
the face numbers. 

The study of these invariants culminated in the proof of the ``$g$-theorem",
which was originally conjectured by McMullen:
\begin{thm}[\cite{BilLee,St80}]\label{g-theorem}
A sequence $f_0, \dots, f_{d-1}\in \Z_{\ge 0}$
are the face numbers of a $d$-dimensional simplicial polytope if and only if 
\begin{itemize}
\item the Dehn-Sommerville relations $h_k = h_{d-k}$ hold,
\item $g_k(P) \ge 0$ for $0 \le k \le d/2$, and
\item there exists a graded ring 
$H = \oplus_{i\ge 0} H_i$, with $H_0 = \R$ and generated by $H_1$, 
for which $g_i = \dim_\R H_i$.  
\end{itemize}
\end{thm}
Sequences satisfying the last condition are known as ``M-sequences".  
Being an M-sequence is equivalent to 
a set of non-linear inequalities. 
For instance, $1, g_1, g_2$ forms an M-sequence if and only if
$g_2 \le \binom{g_1 + 1}{2}$.

The sufficiency of the conditions in Theorem \ref{g-theorem} was
proved by Billera and Lee \cite{BilLee}, who constructed 
appropriate polytopes for every M-sequence.  Necessity 
was first proved by Stanley \cite{St80}, using the fact
that when $P$ is a simplicial polytope, $h_k(P)$ is the 
dimension of the $2k$th cohomology
group of an associated toric variety $X_P$.
This variety is projective, and rationally smooth,
so the Hard Lefschetz theorem implies that $g_k(P)$ is the
dimension of the $2k$th graded piece of the quotient
of $H^\udot(X_P;\R)$ by the ideal generated by 
an element $\lambda\in H^2(X_P;\R)$.  Since the cohomology 
ring of a projective toric variety is generated by the 
elements of degree $2$, so its the quotient by $\lambda$,
which shows that the $g$-numbers are an M-sequence.

Note that although $X_P$ is only defined when $P$ is rational, by
small deformations of the vertices a simplicial polytope can be
made rational without changing its combinatorial type, so Stanley's
argument established the $g$-theorem for all simplicial polytopes.
This deformation trick does not work for general
nonsimplicial polytopes, since there are combinatorial types
of polytopes which cannot be realized over the rationals; see
\cite{Z}.
McMullen later gave a proof \cite{McM} of the necessity part of the
$g$-theorem which did not involve toric varieties and so
worked for non-rational polytopes without a deformation.

\subsection{The toric polynomials}
When $P$ is not simplicial, 
the definition \eqref{simplicial h formula} does
not behave well ---
the resulting polynomial can have negative coefficients,
and the Dehn-Sommerville relations do not hold.  
There is a related problem on the topological
side --- the cohomology of a toric variety defined
by a non-simplicial polytope can exist in odd degrees, 
and the Betti numbers are not invariants of the combinatorics
of faces \cite{McC}. 

The correct way to generalize the simplicial situation 
is to replace the cohomology of $X_P$ with
the intersection (co)homology 
$IH^\udot(X_P;\R)$.  These groups, defined by 
Goresky and MacPherson, are just the usual cohomology
when $X_P$ is rationally smooth, but they are often 
better behaved than cohomology when the variety is singular.
Many important results which hold for cohomology of 
smooth varieties generalize to intersection cohomology of singular varieties.
For instance, the intersection
cohomology Betti numbers of the toric variety $X_P$ depend only 
on the face lattice of $P$, even when $P$ is not simplicial.

The dimensions of these groups are given by the the toric $h$-polynomial,
which Stanley introduced and analyzed in \cite{St87}.  The definition
is recursive: $g(P, t)$ is still computed
from $h(P, t)$ by the formula \eqref{g from h}, but the formula for
$h(P, t)$ becomes 
\begin{equation}\label{h from g}
h(P,t) = \sum_{F<P} g(F,t)(t-1)^{d - 1 - \dim(F)},
\end{equation}
where the sum is over all faces $F$ of $P$, including the empty face
but not $P$ itself.  The induction starts by setting 
$g(P, t) = h(P, t) = 1$ when $P$ is the empty polytope.
It is easy to see that when $P$ is a simplex, $g(P, t) = 1$, so this
definition agrees with the earlier definition for 
simplicial polytopes.

The coefficients of these polynomials depend on more than the face numbers
of $P$; they are $\Z$-linear
combinations of the {\em flag numbers\/} of $P$.  When 
$S\subset \{1,\dots, d-1\}$, the flag number $f_S = f_S(P)$ counts
the number of chains
\[F_{i_1} < \dots < F_{i_r}\]
of faces of $P$, where the elements of $S$ are $i_1 < \dots < i_r$.
When $P$ is simplicial, the flag numbers are determined by the
face numbers, and so are determined by the $g$-numbers.  For general
polytopes, however, there is a lot more information in the
flag numbers than in the $g$-numbers, although the flag
numbers are determined by more general ``convolutions" of the 
$g$-numbers \cite{Kal:new_basis}.

If $P$ is $d$-dimensional, then formulas for the first three 
$g$-numbers are $g_0(P) = 1$ and 
\begin{align}
\label{g1} g_1(P) & = f_0 - (d+1)\\
\label{g2} g_2(P) & = f_1 + f_{02} - 3f_2 - df_0 + \binom{d+1}{2}.
\end{align}
Formulas for higher $g$-numbers in terms of flag-numbers quickly 
become very complicated, but simpler expressions can be obtained by
using a different basis coming from the $\mathbf{cd}$-index \cite{BayEhr}.

Stanley showed \cite{St87} that the toric $h$-numbers satisfy the
Dehn-Sommerville equalities: $h_i(P) = h_{d-i}(P)$ for all $0\le i\le d$.  
For rational polytopes, they follow from Poincar\'e
duality for intersection cohomology, but Stanley gave a 
purely combinatorial proof which works for arbitrary Eulerian posets.
While this gives all linear relations among the face numbers of
simplicial polytopes, for arbitrary polytopes a much larger
set of ``generalized Dehn-Sommerville" relations \cite{BayBil85}
is needed to give all linear relations among flag numbers.

The toric $g$-numbers of $P$ are the dimensions of the graded pieces of
the primitive
intersection cohomology of $X_P$, which is the quotient
of $IH^\udot(X_P; \R)$ by the action of a Lefschetz class in $H^2(X_P; \R)$
(note that intersection cohomology does not have a ring structure, 
but it is a graded module over the cohomology ring).  This
implies the nonnegativity of the $g$-numbers when $P$ is rational. 
The theory of combinatorial intersection cohomology
has extended this to arbitrary polytopes.

Another important inequality among the $g$-numbers is
monotonicity, which was originally
conjectured by Kalai:
\begin{thm}\label{KM}
for any face $F$ of $P$,
\begin{equation*} g(P,t) \ge g(F,t)g(P/F,t),
\end{equation*}
where the inequality is taken coefficient by coefficient.
\end{thm}
Here $P/F$ denotes the ``quotient polytope"
whose face poset is isomorphic to the interval between $F$ and $P$.
This was proved for rational polytopes in \cite{BM} using
a localization argument on intersection cohomology sheaves.
It was pointed out in \cite{BBFK,BrLu} that this argument
translates directly to the combinatorial
IH setting, so the theorem is now known for
arbitrary polytopes.  We give a partial generalization of this
result in Theorem \ref{generalized Kalai} below.

The lack of a ring structure on intersection cohomology
leaves open the question of whether the toric $g$-numbers
give an M-sequence.  Besides the fact that it holds for simplicial
polytopes, there are two other indications that the
answer could be positive.  First, we have the
generalized upper bound theorem:
\begin{thm}[\cite{Bay98}] \label{UBT} For any $d$-polytope $P$ and $1 \le i \le d/2$,
\[g_i(P) \le \binom{f_0 - d + i - 2}{i}.\]
\end{thm}
Bayer's proof in \cite{Bay98} relied on Kalai's monotonicity, which at
the time was only established for rational polytopes.
Now, however, Theorem \ref{KM} is known for arbitrary polytopes, so
Theorem \ref{UBT} is as well.
Theorem \ref{UBT} also follows from a result of Stanley
\cite[Theorem 7.9(a)]{St92}
which says that the $g$-numbers can only increase under subdivisions.
His proof was valid only for rational polytopes, but it
translates immediately to the combinatorial intersection
cohomology setting, so again rationality is not necessary.

The second indication the the toric $g$-numbers might be an $M$-sequence
is the following unpublished result
of Kalai, which is another consequence of Theorem \ref{KM}. 
\begin{thm}  \label{Kalai's theorem}
For a $d$-polytope $P$, $g_k(P) = 0$ implies $g_{k+1}(P) = 0$.
\end{thm}
We present Kalai's ingenious proof in an appendix.

In general we seem to be far from the goal of characterizing all
possible flag $f$-vectors, but many necessary conditions are now known.    
Since the linear equalities are known \cite{BayBil85}, the next 
level of complexity is given by linear inequalities. A summary of the
best linear inequalities currently known for dimensions 
$2$ through $8$ appears in Ehrenborg \cite{Ehr05}.
The inequalities $g_k\ge 0$ and combinations of them
obtained by Kalai's convolution operation \cite{Kal:new_basis}
provide one main source.
Another is the $\mathbf{cd}$-index; we will 
not discuss this here, but note that Karu \cite{Kacd,Kacd2} has recently given a 
construction similar to the combinatorial intersection cohomology
which produces the $\mathbf{cd}$-index.

\section{combinatorial intersection cohomology of fans}

\subsection{Cones and fans}  Although we have described the $h$-polynomial as 
an invariant of polytopes, the combinatorial intersection
cohomology is most naturally described in terms of 
of an associated fan, which carries the same combinatorial structure
but has more convenient geometric properties.   
We first fix some notations and definitions 
regarding convex cones and fans.  

Let $V$ be a finite-dimensional real vector space.
A polyhedral cone in $V$ is a subset of the form
\[\R_{\ge 0}v_1 + \dots+ \R_{\ge 0}v_n,\,\;v_1,\dots,v_n\in V.\] 
All our cones will be assumed
to be pointed, meaning that they do not contain any 
vector subspace other than $\{0\}$.   Any pointed cone $\sigma$
can be expressed as the cone $cP$ over a polytope $P$ of dimension
$\dim P = \dim \sigma - 1$ ---  scale the vectors $v_1,\dots,v_n$
so that they lie in a hyperplane not containing the origin, and let
$P$ be their convex hull.  The resulting polytope is well-defined
up to projective equivalence. 
The map $F \mapsto cF$ is an order-preserving bijection 
between faces 
of $P$ and faces of $\sigma = cP$, where the empty face $\emptyset$ is considered
to be a face of $P$ but not of $\sigma$.  Unless we specifically exclude it,
we always consider a cone or polytope to be a face of itself.

 A {\em fan} $\Delta$ in $V$ is a finite collection of cones in
$V$ so that every face of a cone in $\Delta$ is again in $\Delta$, 
and the intersection of any two cones in $\Delta$ is a face of each.
We use the notation $\tau \prec\sig$ to indicate that a cone
$\tau$ is a face of $\sig$. 
The {\em support} $|\Delta|$ of $\Delta$ is the union
of all its cones.  For instance, 
the set of all faces of a cone $\sig$ forms a
fan with support $\sig$; we denote this fan by $[\sig]$.  

A {\em subfan} of a fan $\Delta$ is a subset which is itself a fan.
For example, given a cone $\sig$, its {\em boundary} 
$\bdy\sig = [\sig]\setminus\{\sig\}$ is a subfan of $[\sig]$.
The {\em relative interior} 
$\sig^\circ$ of $\sig$ is defined to be $\sig \setminus |\bdy\sig|$. 

A fan $\Delta$ in $V$ is {\em complete} if
$|\Delta| = V$.  A polytope $P$ in $V$ gives rise to a complete 
fan $\Delta_P$ in $V$, known as the {\em central fan} of $P$, 
by choosing the origin to be an interior 
point of $P$ and coning off all faces $F \ne P$. 

Given a fan $\Delta$ and $k\ge 0$, its $k$-skeleton $\Delta_{\le k}$ 
is the subfan consisting of all cones of dimension $\le k$, and we write
$\Delta_k$ for the set (not a subfan!) of all cones in $\Delta$ of dimension exactly 
$k$.

If $\Delta$ is a fan in $V$ and $\Delta$ is a fan in $V'$, then
a morphism of fans $\phi\colon \Delta\to \Delta'$ is a linear
map $V\to V'$ so that for every 
cone $\sig \in \Delta$ there exists $\tau \in \Delta'$ with
$\phi(\sig)\subset \tau$.  If $\tau$ can always be chosen so that
$\phi$ induces an isomorphism between $\sig$ and $\tau$, we call
$\phi$ a {\em conewise linear isomorphism}.  One can also define
maps between fans which do not come from a global linear map,
but we will not need them.

\subsection{Conewise polynomial functions for simplicial fans}
For any fan $\Delta$, 
we say a function $f\colon |\Delta| \to \R$
is conewise polynomial if for all $\sig \in \Delta$
the restriction $f|_\sig$ is a polynomial.  The set
$\cA(\Delta)$ of such functions is a graded ring 
under the operations of pointwise addition and multiplication. 
Here we use the usual grading
where linear functions have degree $1$. 

Suppose that $\Delta$ is a simplicial fan, meaning that 
each cone is a cone over a simplex.  As
an abstract poset, $\Delta$ is isomorphic to a simplicial
complex, and the ring $\cA(\Delta)$ is isomorphic to the face ring of this complex,
also known as the Stanley-Reisner ring.   This is a quotient
of the polynomial ring $\R[x_\rho]$ with one generator for
each $1$-cone $\rho\in \Delta_1$ by the ideal 
\[\la x_{\rho_1}\cdots x_{\rho_k} \mid
 \rho_1,\dots,\rho_k\; \text{are not the $1$-faces of a cone in $\Delta$}\ra.\]
The isomorphism comes by identifying
the generator $x_\rho$ with a conewise linear function which 
restricts to a nonzero function on $\rho$ and to zero on
all other $\rho'\in \Delta_1$.  

A simple inclusion-exclusion argument on monomials in the $x_\rho$
shows that the Hilbert series of $\cA(\Delta)$ is
\[\mathop{\rm Hilb}(\cA(\Delta), t) = h(\Delta, t)/(1-t)^{\dim V},\]
where $h(\Delta, t) = \sum_i |\Delta_i|(t-1)^{\dim V - i}$ is the corresponding 
$h$-polynomial.  When $\Delta = \Delta_P$ is the central fan 
of a simplicial polytope $P$, then $h(\Delta, t) = h(P,t)$.  

Realizing the face ring of a simplicial complex as the ring of 
conewise polynomial functions on a fan gives an extra structure
to the face ring in the following way.  Let $A = \Sym(V^*)$ be the ring of 
polynomial functions on $V$.  There is a natural ring homomorphism
$A \to\cA(\Delta)$ obtained by restricting polynomials to $|\Delta|$;
it is injective if $\Delta_{\dim V} \ne \emptyset$.
This makes $\cA(\Delta)$ into an algebra over $A$, and in particular 
a graded $A$-module.

For any graded $A$-module $M$ we define $\ol{M} = M/\m M = M\otimes_A A/\m$,
the quotient by the maximal ideal generated by $A_1$.  
Now assume that $\Delta$ is a $d$-dimensional complete simplicial fan.  
Then one can show that $\cA(\Delta)$ is a free $A$-module.  Since the 
Hilbert series of $A$ is $(1-t)^{-\dim V}$, this gives
\[{\Hilb}(\ol{\cA(\Delta)}, t) = h(\Delta,t).\]

This shows that the $h$-numbers of a complete fan $\Delta$ are nonnegative;
since $h(\Delta_P, t) = h(P, t)$, this gives the nonnegativity of
the $h$-numbers of a polytope $P$.  By the Dehn-Sommerville relations, we have 
\[\dim \ol{\cA(\Delta)}_k = \dim \ol{\cA(\Delta)}_{d-k},\;
0 \le k \le d.\]  Brion \cite{Brion} showed that 
this can be lifted to an $A$-linear dual pairing
\[\cA(\Delta)\otimes_A \cA(\Delta) \to A[d],\]
where $[d]$ shifts the degree down by $d$. 
The pairing is well-defined up to multiplication by a positive
scalar. 

The nonnegativity of the $g$-numbers follows from the existence of a
{\em Lefschetz element} for $\cA(\Delta)$, which is an element
$\ell \in \cA(\Delta)_1$ for which the multiplication 
\[\ell^{d-2k}\cdot\; \colon\ol{\cA(\Delta)}_k \to \ol{\cA(\Delta)}_{d-k}\]
is an isomorphism for $0\le k < d/2$.  This implies that 
$\ell\cdot \colon \ol{\cA(\Delta)}_k \to \ol{\cA(\Delta)}_{k+1}$
is an injection for $k < d/2$ and a surjection for $k > d/2 - 1$, and it follows that 
the quotient ring $H = \ol{\cA(\Delta)}/\ell\,\ol{\cA(\Delta)}$ has
Hilbert series $g(\Delta,t)$.  Thus the existence of a Lefschetz element implies
that the $g$-numbers are nonnegative.  It also shows that they form an M-sequence, 
since $H$ is generated by elements of degree $1$.

When $P$ is a rational simplicial polytope the ring $\ol{\cA(\Delta)}$
is canonically isomorphic to the cohomology ring $H^\udot(X_P; \R)$, by a
map which doubles degree.  The variety $X_P$ is projective, and an
embedding into a projective space determines an ample class in $H^2(X_P;\R)$.  
The Hard Lefschetz theorem (or more precisely its 
extension to rationally smooth varieties due to Saito and 
Beilinson-Bernstein-Deligne) implies that this class is a Lefschetz
element.  This was essentially Stanley's argument to prove the 
necessity part of the $g$-theorem.

McMullen later gave a more elementary proof without using algebraic geometry. 
A conewise linear function  $\ell$ is called 
``strictly convex" if it is convex and it gives a different linear function 
on each full-dimensional cone. Such a function exists if and only if $\Delta$ 
is polytopal, meaning that it is the central fan $\Delta_P$ of a convex polytope $P$.  This is because
the graph of a strictly convex function is the boundary of a cone isomorphic to
the cone $cP$.  
\begin{thm}[McMullen \cite{McM}]\label{simplicial Lefschetz}
If $\Delta$ is a complete simplicial fan, then any strictly convex function
$\ell\in \cA(\Delta)_1$ is a Lefschetz element for $\cA(\Delta)$. 
\end{thm}
Strictly convex linear functions correspond exactly to ample classes in 
$H^2(X_P;\R)$ when $P$ is rational, so McMullen's result gives an elementary
proof of the Hard Lefschetz theorem for rationally smooth toric varieties.

McMullen actually proved the stronger Hodge-Riemann inequalities, which
say that under the Poincar\'e pairing between $\ol{\cA(\Delta)}_k$
and $\ol{\cA(\Delta)}_{d-k}$ the map $(-1)^k\ell^{d-2k}$ is positive definite
on the kernel of $\ell^{d-2k+1}$.  
His proof used induction on dimension, and these
inequalities were essential to make the induction work.   McMullen's
original argument used a different algebraic structure, the
polytope algebra, in place of the ring $\ol{\cA(\Delta)}$; in this
language the Hodge-Riemann inequalities become beautiful statements about
mixed volumes and Minkowski geometry.  The relation between the
polytope algebra and cohomology of toric varieties is described in
\cite{Brion,FulSt}.  Timorin \cite{Tim}
later gave a nice presentation of McMullen's argument
using another description of the ring $\ol{\cA(\Delta)}$ due to
Khovanskii and Pukhlikov.

\begin{rmk}
For the Hodge-Riemann inequalities to hold, it is essential that $\ell$ be strictly
convex.  The condition that $\ell$ be a Lefschetz element is considerably
weaker, however.  In fact, the set of Lefschetz elements is a Zariski open subset of
$\cA(\Delta)_1$, so if it is nonempty, then it is dense.  This suggests that it might
be possible to find Lefschetz elements under weaker assumptions than strict
convexity.  This would be useful to show that 
the conditions of the $g$-theorem 
hold for more general simplicial spheres.
\end{rmk}

\subsection{Example}
Let $\Delta$ be the fan in $\R^2$ whose maximal cones are the four 
quadrants $\{(x,y)\mid \pm x \ge 0, \pm y \ge 0\}$.
Then $h(\Delta, t) = 1 + 2t + t^2$, and the functions 
\[1,\, |x|,\, |y|,\, |xy|\]
give an $A$-module basis for $\cA(\Delta)$.  To see this, 
it is easier to show that
\[1,\, |x|-x,\, |y|-y,\, (|x|-x)(|y|-y)\]
is a basis --- one can subtract off multiples of these
to cancel any element, one quadrant at a time, and this
can be done uniquely.

The strictly convex function $\ell = |x| + |y|$ gives a Lefschetz element; in this case
this just means that $\ell^2 = x^2 + y^2 + 2|xy|$ gives a nonzero element
in $\ol{\cA(\Delta)}_2$.

\subsection{Sheaves on fans}  For a general fan $\Delta$, 
the ring $\cA(\Delta)$ is not as well-behaved as it 
is for simplicial fans. For instance, it is not a free $A$-module, 
and $\ol{\cA(\Delta)}_k$ and
$\ol{\cA(\Delta)}_{d-k}$ may not have the same dimension.
The correct generalization from the simplicial case is the theory
of combinatorial intersection cohomology of \cite{BBFK, BrLu}.  It
is expressed in terms of sheaves on the fan $\Delta$, and we
briefly review how this formalism works.   

We consider a fan $\Delta$ as a topological space
by taking as open subsets all of its subfans.  Given an arbitrary 
subset $\Sig\subset \Delta$, there is a 
unique smallest subfan/open set of $\Delta$ containing $\Sig$, 
namely $[\Sig] = \bigcup_{\sig\in \Sig} [\sig]$.
The closure of a single cone $\{\sig\}$ in this topology 
has a familiar description:
it is the {\em star} 
\[\st(\sig)= \st_\Delta(\sig) = \{\tau\in \Delta\mid \sig \prec \tau\}.\]

A sheaf $\cF$ of vector spaces 
on $\Delta$ is given by an assignment $\Sig \mapsto \cF(\Sig)$ of a 
vector space to each
subfan $\Sig\subset \Delta$, together with restriction maps
$\cF(\Sig)\to \cF(\Sig')$ for any pair of subfans $\Sig'\subset \Sig$.
These are required to satisfy
(1) if $\Sig'' \subset \Sig' \subset \Sig$ are subfans of
$\Delta$, then the obvious triangle of restrictions is commutative,
and (2) if $\Sig, \Sig'$ are any two subfans of $\Delta$, then
the image of $\cF(\Sig\cup\Sig')$ in $\cF(\Sig)\oplus\cF(\Sig')$
is the set of pairs $(x,y)$ for which $x$ and $y$ restrict to 
the same element of $\cF(\Sig\cap\Sig')$.  Sheaves of rings
or modules are defined in the same way.

For any fan $\Delta$, the assignment $\Sig \mapsto \cA(\Sig)$
gives a sheaf known as the sheaf of conewise 
polynomial functions; we denote it by $\cA_\Delta$ or simply
$\cA$ if the fan is understood.  It is a sheaf of rings, which 
means we can define sheaves of modules over it, as follows.  
An $\cA_\Delta$-module is a sheaf $\cF$ on 
$\Delta$ together with a structure of a graded $\cA(\Sig)$-module
on $\cF(\Sig)$ for every subfan $\Sig\subset \Delta$.
This is required to be compatible with the restriction maps:
$\cF(\Sig')\to\cF(\Sig)$ should be a 
homomorphism of $\cA(\Sig')$-modules whenever 
$\Sig\subset\Sig'$ are subfans of $\Delta$, where $\cF(\Sig)$
becomes a $\cA(\Sig')$-module via the restriction
$\cA(\Sig') \to \cA(\Sig)$. 

Because $\Delta$ is finite, a sheaf $\cF$ on $\Delta$
can be described by a finite collection of modules and maps.  
The stalk of $\cF$ at a cone $\sig\in \Delta$ is 
$\cF([\sig])$, since $[\sig]$ 
is the smallest open set containing $\sig$.  To simplify notation, we will 
write $\cF(\sig)$ for $\cF([\sig])$.  If $\tau\prec \sig$,
then restriction of sections gives a homomorphism $\cF(\sig)\to \cF(\tau)$.
Since any subfan of $\Delta$ is a union of fans of the form $[\sig]$,
the data in the sheaf $\cF$ is equivalent to 
the collection of stalks together with the restriction maps 
between them (subject to the obvious commutation relation). 

The stalks of $\cA_\Delta$ are particularly simple: $\cA_\Delta(\sig)$ is the ring  
of polynomial functions on the span of $\sig$, which we denote
by $A_\sig$.  If $\cF$ is an $\cA_\Delta$-module,
then the $\cF(\sig)$ is a $A_\sig$-module for every face $\sig$.

For a $\Delta$-sheaf $\cF$ and a pair of subfans $\Sig\subset \Sig'$ of $\Delta$
we define the space of relative sections $\cF(\Sig',\Sig)$ to be the kernel
of the restriction $\cF(\Sig')\to \cF(\Sig)$.

If $\phi\colon \Delta\to \Delta'$ is a morphism of fans, and
$\cF$ is a sheaf on $\Delta$, then the pushforward sheaf
$\phi_*\cF$ on $\Delta'$ is defined by 
$\phi_*\cF(\Sig) = \cF(\phi^{-1}(\Sig))$, where
\[\phi^{-1}(\Sig) = \{\tau \in \Delta \mid \phi(\tau) \subset |\Sig|\}.\]
If $\phi$ is a conewise linear isomorphism, then 
$\phi_*\cA_{\Delta} \cong \cA_{\Delta'}$; in general we
only have a natural map $\cA_{\Delta'} \to \phi_*\cA_{\Delta}$.
If $\cF$ is an $\cA_{\Delta}$-module, then this map makes $\phi_*\cF$ 
into an $\cA_{\Delta'}$-module.
\subsection{Combinatorial IH sheaves}  
The main definition in the theory of combinatorial intersection cohomology 
developed in \cite{BBFK, BrLu} is the following.

\begin{defn}
An $\cA_\Delta$-module $\cF$ is called {\em pure} if it is 
\begin{itemize}
\item {\em locally free}, meaning that $\cF(\sig)$ is a free 
$A_\sig$-module for every $\sig\in\Delta$, and
\item {\em flabby}, meaning that the restriction
$\cF(\Delta)\to \cF(\Sig)$ is surjective for any subfan $\Sig$ of $\Delta$.
\end{itemize}
\end{defn}
Note that flabbiness of $\cF$ is equivalent to either
\begin{enumerate}
\item $\cF(\Sig)\to \cF(\Sig')$ is surjective for any subfans $\Sig'\subset \Sig$ of $\Delta$, or\
\item $\cF(\sig)\to \cF(\bdy\sig)$ is surjective for every $\sig\in \Delta$. 
\end{enumerate}

We can construct indecomposable pure sheaves inductively as follows.  
Starting with a cone $\sig\in\Delta$,
define a sheaf $\cL = {}_\sig\cL_\Delta$ to be zero on all 
cones $\tau\notin \st(\sig)$, and let $\cL(\sig) = A_\sig$.
Then, assuming that $\cL$ has already been defined on $\bdy\tau$
for $\tau\in \Delta$, let $\cL(\tau)$ be a minimal free $A_\tau$-module which 
surjects onto $\cL(\bdy\tau)$, where $\cL(\bdy\tau)$ is an $A_\tau$-module 
via the restriction $A_\tau = \cA(\tau) \to \cA(\bdy\tau)$. 
Equivalently, let $\cL(\tau)$ be a free $A_\tau$-module generated
by a vector space basis for $\ol{\cL(\bdy\tau)}$, and let the map
$\cL(\tau) \to \cL(\bdy\tau)$ be given by choosing representatives
for the basis elements.   

The following ``decomposition theorem" is an easy consequence of these
definitions.
\begin{thm}[\cite{BBFK, BrLu}] \label{semisimplicity}
The sheaves ${}_\sig\cL_\Delta$, $\sig\in \Delta$ give 
a complete list of isomorphism classes of indecomposable
pure $\cA_\Delta$-modules, up to shifts of degree.
Every pure sheaf is isomorphic to a direct sum of 
indecomposable objects.
\end{thm}
The sheaves ${}_\sig\cL_\Delta$ are the combinatorial IH sheaves.
For many questions it is enough to look at ${}_\zcone\cL_\Delta$, where
$\zcone = \{0\}$ is the zero cone.  This is because
the other combinatorial IH sheaves ${}_\sig\cL_\Delta$ 
can be derived from it by a pullback and base change
from a quotient fan.  To simplify notation, 
we put $\cL_\Delta = {}_\zcone\cL_\Delta$.
Also, note that if $\Sig\subset\Delta$ is a subfan containing the cone 
$\sig$, then  the restriction of ${}_\sig\cL_\Delta$ to $\Sig$ is 
${}_\sig\cL_\Sig$, so we can write ${}_\sig\cL(\Sig)$ for the
sections on $\Sig$ without causing confusion.

If $\Delta$ is simplicial, then $\cA_\Delta$ is pure,
so $\cL_\Delta$ is (canonically) isomorphic to
$\cA_\Delta$.  Thus global
sections of $\cL$ are just conewise polynomial functions;
if $\Delta = \Delta_P$ for a simplicial polytope $P$, then
$\Hilb(\ol{\cL(\Delta)}, t)$ is the simplicial $h$-polynomial
defined by \eqref{simplicial h formula}.   

This generalizes to arbitrary fans, using the toric $g$ and $h$-numbers:
\begin{thm}\label{big theorem} Let $P$ be a polytope.
\begin{itemize}
\item[(a)] If $\Delta = \Delta_P$ is the central fan of a polytope $P$,
and $\sig = cF$ is the cone over a proper face $F$ of $P$, 
then ${}_\sig\cL(\Delta)$ is a free $A$-module and
\[\Hilb(\ol{{}_\sig\cL(\Delta)},t) = h(P/F,t).\]
\item[(b)] If $\tau = cP$ and $\sig = cF$ for a face
$F \le P$, then 
$\Hilb(\ol{{}_\sig\cL(\tau)}, t) = g(P/F,t)$.
\end{itemize}
\end{thm}
Unlike the simplicial case, the proof of this is far from
straightforward.  In order to explain the basic structure of
the proof, we restrict for simplicity to the case
$\sig=\zcone$.  Let (a)${}^{}_d$
and (b)${}^{}_d$ represent the statements (a) and (b), 
restricted to polytopes $P$ of dimension $\le d$.  Theorem \ref{big theorem} 
is proved by a spiraling induction, 
by showing that (b)${}_{d-1} \implies$ (a)${}^{}_d$ and 
(a)${}^{}_d \implies$ (b)${}^{}_d$.

The first implication follows from the statement that for a complete
fan the cellular cohomology gives an exact sequence
\[0 \to \cL(\Delta) \to \bigoplus_{\sig \in \Delta_d} \cL(\sig) \to \bigoplus_{\tau \in \Delta_{d-1}} \cL(\tau) 
\to \dots \to \cL(\zcone) \to 0,\]
where the maps are sums of all possible restriction maps, with appropriate
minus signs added to make it a complex.
Together with some simple commutative
algebra this implies that $\cL(\Delta)$ is a free $A$-module.
Taking Hilbert series, we obtain
\[(1-t)^{-d}\Hilb(\ol{\cL(\Delta_P)},t) = 
\sum_{\sig \in \Delta} (-1)^{d - \dim \sig}\Hilb(\ol{\cL(\sig)},t)(1-t)^{-\dim \sig},
\]
which shows that $\Hilb(\ol{\cL(\Delta_P)},t) = h(P, t)$, using  
\eqref{h from g} and the assumption (b)${}_{d-1}$.

The other implication (a)${}^{}_d \implies$ (b)${}^{}_d$ is more subtle.
Take $P$ a $d$-polytope, and 
let $\sig = cP$.  Let $V$ be the linear span of $\sig$.  
If $v$ is in the relative interior of
$\sig$, then the projection $\pi\colon V \to W = V/\R v$ gives a
conewise linear isomorphism from $\bdy\sig$ to a complete 
fan $\Delta$ in $W$ which is isomorphic to 
the central fan $\Delta_P$.  Conversely, 
$\bdy\sig$ can be viewed as the graph of a strictly convex conewise linear function 
$\ell \in\cA(\Delta)_1$.

Since $\pi$ induces a linear isomorphism between cones of $\bdy\sig$
and cones of $\Delta$, there is an isomorphism of $\cA_{\Delta}$-sheaves
$\pi_*\cL_{\bdy\sig} \cong \cL_\Delta$.
Taking global sections, we get an isomorphism
\begin{equation}\label{projection iso}
\cL(\bdy\sig) \cong \cL(\Delta)
\end{equation}
of $A'$-modules, 
where $A' = \Sym(W^*)$ is included as a subring of $A = \Sym(V^*)$ 
by the pullback $\pi^*$.  
Choosing a degree $1$ element 
$y\in A \setminus A'$ gives an isomorphism $A \cong A'[y]$.  In
terms of the isomorphism \eqref{projection iso}, the action 
of $y$ on $\cL(\Delta)$  is 
given by multiplication by $\ell$. 

By the construction of $\cL$, the map
$\ol{\cL(\sig)} \to \ol{\cL(\bdy\sig)}$ is an isomorphism.
On the other hand, the isomorphism \eqref{projection iso} 
implies that $\ol{\cL(\bdy\sig)} \cong \ol{\cL(\Delta)}/\ell\,\ol{\cL(\Delta)}$. 
The inductive hypothesis (a)${}_d$ says that
the Hilbert series of $\ol{\cL(\Delta)}$ is $h(P,t)$.  
Thus (b)${}_d$ follows from the following difficult result of Karu.
\begin{thm}[Karu \cite{Ka}]\label{Karu}
If $\Delta=\Delta_P$ for a polytope $P$, then any strictly convex 
conewise linear function $\ell$ acts on $\ol{\cL(\Delta)}$ 
as a Lefschetz operator.
\end{thm}
As with McMullen's Theorem \ref{simplicial Lefschetz},
what Karu actually proved is the stronger 
Hodge-Riemann inequalities; these are needed
to make the induction work.  The inequalities are taken
with respect to a dual pairing 
\[\ol{\cL(\Delta)}_k \otimes \ol{\cL(\Delta)}_{d-k} \to \R.\]
Karu used a pairing from \cite{BBFK} which involved choices, 
but Bressler and Lunts \cite{BrLu05} have simplified Karu's argument 
by using a canonical pairing they defined
in \cite{BrLu}.  Another presentation of Karu's theorem
appears in \cite{BBFK3}, using 
another approach to defining a canonical
pairing from \cite{BBFK2}.  

Karu's theorem implies the following degree vanishing
result.
\begin{thm} \label{degree vanishing}
For $\sig \prec \tau$, the stalk ${}_\sig\cL(\tau)$
is generated in degrees $< (\dim \tau - \dim \sig)/2$,
while the costalk \[{}_\sig\cL{(\tau,\bdy\tau)}
= \ker({}_\sig\cL(\tau) \to {}_\sig\cL(\bdy\tau))\] is 
generated in degrees $> (\dim\tau-\dim\sig)/2$.
\end{thm}

In particular, ${}_\sig\cL(\tau) \to {}_\sig\cL(\bdy\tau)$
is an isomorphism in degrees $< (\dim \tau - \dim \sig)/2$.
This has the following important consequence.
\begin{cor}\label{scalar automorphisms}
 The only automorphisms of the combinatorial IH sheaf ${}_\sig\cL$ (as a
graded $\cA_\Delta$-module) are multiplication by scalars $\R^\times$. 
\end{cor}
This means that the combinatorial IH groups are 
{\em canonically} associated to the fan.  
 As a result, the groups
themselves carry interesting information about $P$ 
beyond just their dimensions.  

\subsection{Connections with topology}
Although the theory of combinatorial IH sheaves can be 
developed without referring to toric varieties, the
topological interpretation still provides a powerful
way of understanding what these results mean.
All the cohomology and intersection cohomology spaces
below are taken with $\R$ coefficients.

Suppose that the fan $\Delta$ is rational with 
respect to a fixed lattice in the vector space $V$.
The definition of the combinatorial IH sheaf $\cL_\Delta$ 
was first given by Barthel, Brasselet, Fieseler and Kaup
in \cite{BBFK1}, where they show that the $\cL(\Delta)$ is isomorphic to 
the equivariant intersection cohomology 
$IH^{\udot}_T(X_\Delta)$ of the associated toric variety under the 
action of the usual torus $T \cong (\C^*)^d$.

Equivariant intersection cohomology, which was defined in 
\cite{BerLu,Bry,J}, is a topological invariant associated to 
a space endowed with an action of a group.  When the group
is trivial, it specializes to the usual intersection
cohomology, and when the space is smooth, it gives the 
equivariant cohomology.    
The extra information in the group action gives the
equivariant cohomology and intersection cohomology
more structure and better properties than their
non-equivariant counterparts.  In particular, 
they are modules over the equivariant cohomology
of a point $H^\udot_T(p)$, which is canonically isomorphic
to $A = \Sym(V^*)$ by an isomorphism which 
identifies $V^*$ with $H^2_T(p)$ (note that 
our choice to use the standard grading on $A$ means
that all our isomorphisms between the algebraic and 
topological sides will double degree).

If $\Sig\subset\Delta$ is a subfan, then $X_\Sig$ is
an open $T$-invariant subvariety of $X_\Delta$.  
Equivariant IH restricts via open inclusions, so there is
a homomorphism 
$IH^\udot_T(X_\Delta) \to IH^\udot_T(X_\Sig)$ of $A$-modules.  
Barthel, Brasselet, Fieseler and Kaup showed that with these 
restrictions, the assignment 
\begin{equation}\label{IH sheaf}
\Sig\mapsto IH^\udot_T(X_\Sig), \;\Sig\subset \Delta 
\end{equation}
is a sheaf on $\Delta$.

The fact that this is a presheaf is purely formal ---
the same is true for the ordinary cohomology and 
intersection cohomology, for instance. 
But being a sheaf is special to equivariant intersection 
cohomology; for most types of cohomology it is not possible to 
describe classes locally by gluing them from classes
on sets in an open cover. 

The equivariant cohomology also restricts along open inclusions,
but it does not give a sheaf in general
(see \cite{BBFK1} for an example).  However, it
has a sheafification, which is just the sheaf
 $\cA_\Delta$ of conewise polynomial functions.  
In particular, if $\tau$ is a cone, then 
$H_T^\udot(X_{[\tau]}) \cong A_\tau$
canonically.
Equivariant intersection cohomology is always a module
over equivariant cohomology, which makes the sheaf
\eqref{IH sheaf} into an $\cA_\Delta$-module.

To show that this sheaf is the combinatorial IH sheaf
$\cL_\Delta$, one needs to show that $IH^\udot_T(X_{[\tau]})$ is 
the minimal free $A_\tau$-module
which surjects onto $IH^\udot_T(X_{\bdy\tau})$
for every cone $\tau\ne \zcone$.  This 
was shown in \cite{BBFK1}.  It also follows from the
following result of Bernstein and Lunts.
\begin{thm} \cite{BerLu}
Let $X$ be a variety with a $T$-action, and suppose that
there is a homomorphism $\rho\colon \C^*\to T$ which 
contracts $X$ to a point $y$:
\[\lim_{t\to 0} \rho(t)\cdot x = y\;\text{for all $x\in X$}.\]
Then $IH^\udot_T(X)$ is the minimal free $H_T^\udot(p)$-module which surjects onto
$IH^\udot_T(X \setminus \{p\})$.
\end{thm}
The proof in \cite{BerLu} used the Hard Lefschetz theorem for intersection
cohomology.  A different proof, using the weight filtration on equivariant
intersection cohomology, appears in \cite{BM2}.

When $\Delta$ is a complete fan or $\Delta = [\tau]$, the
variety $X_\Delta$ is ``formal" for equivariant intersection
cohomology in the sense of \cite{GKM}.  This means that $IH^\udot_T(X_\Delta)$
is a free $A$-module and
\[IH^\udot(X_\Delta) = \ol{IH^\udot_T(X_\Delta)}\]
canonically.  Barthel, Brasselet, Fieseler and Kaup
\cite{BBFK} gave a complete characterization of when
$X_\Delta$ is formal; we will use it in Section
\ref{positive local h} below.

When $X_\Delta$ is formal, then $\ol{\cL(\Delta)}$ is 
{\em canonically} isomorphic to $IH^\udot(X_\Delta)$.  So 
Karu's Theorem \ref{Karu} is not simply a substitute for
the Hard Lefschetz theorem; it is actually an alternate
proof of it, valid for toric varieties.  
A careful study of Karu's proof and the combinatorial 
IH sheaves in general should give new insight into the Hard Lefschetz
theorem, one of the deepest theorems in algebraic geometry.

This will hopefully lead to similar elementary proofs in other interesting 
contexts.  For instance, \cite{BM2} defines a similar theory of
combinatorial IH sheaves for Schubert varieties; the role of fans
is played by ``moment graphs,'' which are linearly embedded graphs encoding
the fixed points and invariant curves for a torus action.  
At present these results are only valid when the graph actually
comes from a variety, as there is no result comparable to Karu's theorem
so far.  In fact, we don't even have an elementary proof of Hard Lefschetz
for ``smooth" moment graphs, which have been extensively studied 
by Guillemin and Zara \cite{GZ1,GZ2,GZ3,GZ4}.  

There is a 
natural class of non-rational moment graphs for which the 
construction of \cite{BM2} should give good answers: the 
``Bruhat graphs" arising from a non-crystallographic 
Coxeter group $W$.  These can be thought of as the 
moment graphs of the (nonexistent) Schubert varieties in the 
flag variety of the 
(nonexistent) semi-simple group with Weyl group $W$,  
just as non-rational fans correspond to nonexistent toric varieties.
The stalks of the combinatorial IH sheaf on a Bruhat graph for $W$
should be free, with ranks given by Kazhdan-Lusztig polynomials for $W$; 
this would give a proof that these polynomials have nonnegative coefficients.
See \cite{Fieb} for an exposition of these ideas.

Not only does the space of sections of the sheaf $\cL$ give the 
module  $IH^\udot_T(X_\Delta)$ when $\Delta$ is rational, but the sheaf 
$\cL_\Delta$ itself
can be seen as a model for the equivariant IH sheaf $IC^\udot(X_\Delta)$, an
object in the equivariant derived category $D^b_T(X_\Delta)$ \cite{BerLu},
and complexes of pure sheaves can be used to model more general 
objects of $D^b_T(X_\Delta)$.  This will be discussed in more detail in
section \ref{Koszul}.

\subsection{Examples in low degree: degree zero and one}
\label{degree zero and one}
To illustrate the features of this theory, we will describe
the graded pieces of the sheaf $\cL = {}_{\zcone}\cL_\Delta$ in 
degrees $\le 2$.
The answers are given as cohomology of simple chain complexes whose 
terms have dimensions corresponding directly to terms in 
the formulas \eqref{g1} and \eqref{g2} for $g_1$ and $g_2$.

The degree zero part of $\cL$ is trivial; we have
\begin{equation}\label{L in degree 0}
\cL(\Delta)_0 = \R
\end{equation}
for any fan $\Delta$.

In order to understand the degree 
$k$ part of $\cL$ for $k\ge 1$, it is enough
to understand its restriction to the $(2k-1)$-skeleton
$\Delta_{\le 2k-1}\subset \Delta$.  This is because 
Theorem \ref{degree vanishing} implies that the restriction
\[\cL(\Delta)_{k} \to \cL(\Delta_{\le 2k-1})_k\]
is an isomorphism.

For $k=1$, this means that it is enough to consider sections
of $\cL$ on $\Delta_{\le 1}$.  Since this is a simplicial fan,
we have 
\begin{equation}\label{L in degree 1} 
\cL({\Delta})_1 = \cL({\Delta_{\le 1}})_1 = \cA({\Delta_{\le 1}})_1 = 
\bigoplus_{\rho \in \Delta_1} \la \rho\ra,
\end{equation}
where we introduce the notation $\la \tau\ra = (\Span \tau)^*$
for the space of linear functions on a cone $\tau$.
 
Reducing modulo $\m$, we see that $\ol{\cL(\Delta)}_1$ is the cokernel of the
multiplication map $\mu_1\colon V^*\otimes \cL(\Delta)_0 \to \cL(\Delta)_1$. 
In terms of the identifications \eqref{L in degree 0} and \eqref{L in degree 1}, 
this is the map 
\begin{equation}\label{mult 1}
V^* \to \bigoplus_{\rho \in \Delta_1} \la \rho\ra
\end{equation}  
induced from the natural restrictions.
  
If the $1$-cones $\rho \in \Delta_1$ span $V$, 
for instance if $\Delta$ contains a full-dimensional cone, 
then $\mu_1$ is injective, so
$\dim \ol{\cL(\Delta)}_1 = \#\Delta_1 - \dim(V)$. 
Thus for a $d$-polytope $P$ we have
\begin{align*}
\dim \ol{\cL(\Delta_P)}_1 & = f_0(P) - d = h_1(P)\\
\dim \ol{\cL(cP)}_1 & = f_0(P) - (d + 1) = g_1(P),
\end{align*}
as predicted by Theorem \ref{big theorem}.

This implies that $g_1$ and $h_1$ are nonnegative, which is just
the obvious statement that a $d$-polytope must have at least
$d+1$ vertices.  Still, it is interesting to have 
canonical geometrically defined vector
spaces of dimensions $g_1(P)$ and $h_1(P)$.  For instance, 
the dual of $\ol{\cL(cP)}_1$ is canonically
the space of all affine dependencies among the vertices of $P$.  
This interpretation of $g_1$ was used by Kalai to prove a number of
results, including low-degree cases of his monotonicity conjecture.

Note also that the fact that nothing new needs to be added in degree
one for larger cones can be seen directly without using Theorem 
\ref{degree vanishing}.  What is needed is to see that 
for any cone $\tau$, the $A_\tau$-module defined by \eqref{L in degree 0}
and \eqref{L in degree 1} with module structure given by
\eqref{mult 1} is isomorphic to the degree zero and one part
of a free module.  But this is just the obvious fact that
\[\la \tau\ra \to \bigoplus_{\rho\in [\tau]_1} \la \rho\ra\]
is injective for all cones $\tau$.

\subsection{Degree two} \label{degree two} 
To understand $\cL$ in degree two it is enough
to understand its restriction to the $3$-skeleton $\Delta_{\le 3}$.
Since any two-cone is simplicial, $\cL$ will agree with $\cA$ on 
$\Delta_{\le 2}$.  There will be a correction for nonsimplicial $3$-cones,
however, since $\cA$ will not be flabby on such cones.

To describe this correction, fix a choice of a subdivision $\wt\Delta$ of 
$\Delta_{\le 3}$ which divides each 
$3$-cone into simplicial cones without adding new $1$-cones; in particular the
$2$-cones are not subdivided.  Let $\phi\colon \wt\Delta\to \Delta$
be the associated map of fans.  Then the pushforward 
$\phi_*\cA_{\wt\Delta}$ is isomorphic to the restriction of $\cL$ to $\Delta_{\le 3}$.
The restrictions of these sheaves to the two-skeleton $\Delta_{\le 2}$ are
clearly isomorphic, since $\cA_{\Delta_{\le 2}} = \cL_{\Delta_{\le 2}}$
because $\Delta_{\le 2}$ is simplicial.  
Constructing the isomorphism on a 
$3$-cone $\sig$ amounts to showing that $\cA(\wt{[\sig]})$ is a free $A_\sig$-module
generated in degrees zero and one, and that the restriction 
$\cA(\wt{[\sig]}) \to \cA(\bdy\sig)$ is surjective, with 
kernel generated in degrees $\ge 2$.   

By Corollary \ref{scalar automorphisms}, this isomorphism is canonical. 
As a consequence, the sheaves obtained by different choices of
subdivisions $\wt\Delta$ are canonically isomorphic.  This
somewhat surprising fact can be understood as follows.  Let $\cF = \phi_*\cA_{\wt\Delta}$.
If $\sig \in \Delta$ is
a $3$-cone, then the multiplication map
\[V^*\otimes \cF(\sig)_1 \to \cF(\sig)_2\]
is surjective.  This can be checked directly, or else deduced from
Theorem \ref{degree vanishing}, which implies that $\ol{\cF(\sig)}_2 = 0$.
Thus classes in $\cF(\sig)_2$ can be specified by choosing a 
preimage in $V^*\otimes \cF(\sig)_1$, which is clearly independent of
the chosen subdivision. 

Thus $\cL(\Delta)_2$ is identified 
with degree two conewise polynomial functions on the simplicial fan 
$\wt{\Delta}_{\le 3}$.  As we noted earlier, sections of $\cA$ on
a simplicial fan gives the face ring of
the corresponding abstract simplicial complex.  The 
decomposition of $\cL(\Delta)_1$ by the monomial generators
is just \eqref{L in degree 1}; 
the corresponding decomposition of $\cL(\Delta)_2$ by monomials gives
an identification
\begin{equation}\label{L in degree 2}
\cL(\Delta)_2 = \bigoplus_{\rho_1\triangleleft\rho_2} \la\rho_1\ra \otimes \la \rho_2\ra,
\end{equation}
where the relation $\rho_1 \triangleleft \rho_2$ for $\rho_1, \rho_2\in \Delta_1$ 
means that either $\rho_1 = \rho_2$ or 
else $\rho_1$ and $\rho_2$ generate a two-dimensional cone of
$\wt{\Delta}$ and $\rho_1$ precedes
$\rho_2$ in some fixed total order on $\Delta_1$.

The quotient space $\ol{\cL(\Delta)}_2$ is the cokernel of the
multiplication map $\mu_2\colon V^* \otimes \cL(\Delta)_1 \to \cL(\Delta)_2$.
In terms of \eqref{L in degree 1} and \eqref{L in degree 2}, this map is
\[\bigoplus_{\rho\in \Delta_1} V^* \otimes \la \rho\ra \to 
\bigoplus_{\rho_1\triangleleft\rho_2} \la\rho_1\ra \otimes \la \rho_2\ra,\]
where $V^*\otimes \la \rho\ra$ maps to all terms on the right of the form
$\la\rho\ra \otimes \la \rho'\ra$ or $\la\rho'\ra \otimes \la \rho\ra$ via
the natural map $V^*\to \la\rho'\ra$.

Unlike the discussion of $\mu_1$ in the previous section, $\mu_2$ is never injective, because 
$\cL(\Delta)$ has a generator in degree zero.  Writing down the ``obvious"
elements in the kernel of $\mu_2$ gives rise to a two-step chain complex 
\begin{equation}\label{chain complex}
\bigwedge\nolimits^2 V^* \stackrel{\phi}\lra \bigoplus_{\rho\in \Delta_1} V^* 
\otimes \la \rho\ra \stackrel{\mu^{}_2}\lra \bigoplus_{\rho_1\triangleleft\rho_2} 
\la\rho_1\ra \otimes \la \rho_2\ra,
\end{equation}
where \[\phi(\alpha\wedge\beta) = \alpha\otimes \mu_1(\beta) - \beta\otimes\mu_1(\alpha)\]
for any $\alpha, \beta\in V^*$.  In other words, $\mu_2\circ\phi = 0$.

Now suppose that $\Delta$ is either the
central fan $\Delta_P$ 
or the cone $[cP]$ for some polytope $P$.
over a polytope of dimension $d-1$.  
Then the complex
\eqref{chain complex} is left exact; i.e., it has cohomology only at the rightmost
place.  The injectivity of $\phi$ follows from the fact that $\Delta$ contains
a full-dimensional cone, since then $\phi$ factors through the natural map
$\bigwedge^2 V^* \to V^* \otimes V^*$.  Exactness in the 
middle follows from the fact that $\cL(\Delta)$ is a free $A$-module
in positive degrees with a single generator in degree zero. 
Note that this depends on the identifications \eqref{L in degree 1} and
\eqref{L in degree 2}, which in turn required the use of Theorem
\ref{degree vanishing}.

\subsection{Connection with infinitesimal rigidity}
Conversely, knowing the exactness of \eqref{chain complex}
when $\Delta = [\sig]$ for $\sig$ a cone of dimension $d > 3$
implies, without appealing to Theorem
\ref{degree vanishing}, that the restriction 
$\cL(\sig)_2 \to \cL(\bdy\sig)_2$
is an isomorphism, and thus by induction that 
$\cL(\sig)_2 \to \cL(\Delta_{\le 3})_2$ is an
isomorphism, which in turn implies the formula \eqref{L in degree 2}. 
As we will see, this exactness can be deduced
from results of Whiteley about infinitesimal rigidity of frameworks.
So the theory up to degree two can be verified independently,
without using Karu's theorem.

To see the rigidity interpretation of \eqref{chain complex}, we 
replace it by a smaller but equivalent
subcomplex.  Both the source and target of $\mu_2$ surject onto
$\cL(\Delta_{\le 1})_2 = \bigoplus_{\rho \in \Delta_1} 
\la \rho\ra \otimes \la \rho\ra$ via maps commuting with $\mu_2$. 
Taking kernels gives a smaller complex with the same cohomology:
\begin{equation}\label{rigidity complex}
0\to \bigwedge\nolimits^2 V^* \to \bigoplus_{\rho\in \Delta_1} 
\rho^\bot \otimes \la \rho\ra \stackrel{\tilde\mu^{}_2}\lra
\mathop{\bigoplus_{\rho_1\triangleleft\rho_2}}_{\rho_1\ne\rho_2}
\la\rho_1\ra \otimes \la \rho_2\ra,
\end{equation}
where $\rho^\bot \subset V^*$ is the annihilator of $\Span(\rho)$.

Now suppose that $\Delta = [cP]$ is the cone over a $d$-polytope
$P$ in a vector space $W$, so that $\Delta$ is a fan in 
$V = W \oplus \R$.  The dimensions of the 
terms of \eqref{rigidity complex} 
of $P$ are $\binom{d+1}{2}$, $df_0$, and $f_1 + (f_{02} - 3f_2)$,
respectively.  Thus the formula 
\eqref{g2} for $g_2(P)$ computes the Euler characteristic
of this complex.  Since the complex has cohomology only
at the rightmost term, we see that $g_2(P)$ is the 
dimension of $\coker(\tilde\mu_2) \cong \coker(\mu_2) = 
\ol{\cL(\Delta)}_2$, as asserted by Theorem \ref{big theorem}.  
 
The fan $\Delta_{\le 2}$ is the 
cone over a framework, or geometric graph, whose vertices
are the vertices of $P$ and whose edges are the edges of
$P$ together with enough extra edges to triangulate each
$2$-face.
A $1$-cone $\rho\in \Delta_1$ will be the cone over a 
vector $(v,1)$, where $v\in V$ is a vertex of $P$.  
This gives an identification $\la \rho\ra \cong \R$ via 
$\phi \mapsto \phi(v,1)$.  By restricting covectors
we get an identification of $\rho^\bot \subset (W\oplus \R)^*$ 
with $W^*$, and choosing
an inner product on $W$, we get an identification
$\rho^\bot \cong W$.  

With these identifications the map $\tilde\mu_2$ on the right 
of \eqref{rigidity complex} becomes the infinitesimal 
rigidity matrix of the framework, which is a map
\begin{equation}\label{rigidity matrix}
\bigoplus_{v} W \to \bigoplus_{e} \R,
\end{equation}
where the sum is over all vertices $v$ and edges $e$ of the
framework.  It sends a tuple $(w_v)$ to 
$(w_v - w_{v'})\cdot (v - v')$
on an edge $e$ with endpoints $v, v'$.  
Elements of $\bigoplus_{v} W$ are assignments of vectors
to each vertex of $P$ which should be thought of as infinitesimal 
motions of these points.  Being in the kernel of $\mu_2$ means that
the lengths of the edges of $P$ do not change to first order under this
motion, with respect to the chosen inner product.

The image of $\bigwedge^2 V^*$ in $\bigoplus_{v} W$ gives
assignments of vectors which come from global affine motions of $W$.  
Thus the exactness of \eqref{chain complex} translates to the statement
that the framework is infinitesimally rigid --- all motions of the
vertices preserving the lengths of the edges to first order come from
global affine motions.  When $\dim P = 3$, 
this was proved by Aleksandrov, generalizing theorems of 
Cauchy and Dehn for the case when $P$ is simplicial (see also \cite{FP,P}).  
In higher dimensions,
Whiteley \cite{Wh} showed how to deduce rigidity inductively from
the three-dimensional case.  

The cohomology of the dual complex to \eqref{rigidity complex} is
known as the space of stresses of the framework; these can be viewed
as assignments of expanding or contracting forces along the
edges in such a way that the total forces at each vertex cancel.
Before the advent of combinatorial intersection cohomology,
Kalai \cite{Kal:rigidity,Kal:aspects} pointed out that this
gives a vector space of dimension $g_2(P)$, thus proving 
$g_2(P) \ge 0$.  He also used the interpretation of $g_2$ as
a space of stresses to prove the degree two piece of his monotonicity
conjecture: $g_2(P) \ge g_2(F) + g_1(F)g_1(P/F) + g_2(P/F)$.  

Several authors \cite{Lee, Lee2, TWW, TWW2} have 
considered a theory of ``higher rigidity" for simplicial complexes
which plays the same role for higher $g_k$'s that infinitesimal
rigidity plays for $g_2$.  These theories are closely related to 
McMullen's polytope algebra.

Note that the connection between $g_2$ and rigidity explains why
Theorems \ref{Karu} and \ref{degree vanishing} depend on
the inductive nature of fans (they are made up of
cones over convex polytopes, which can be described by 
polytopal fans in one dimension less, and so on). 
Connelly \cite{Conn1,Conn2} has
constructed non-convex simplicial triangulated spheres in $\R^3$ 
which are not rigid.   The cone over such a sphere will be a
$3$-dimensional fan in $\R^4$; if this is considered to be 
the boundary of a non-convex ``cone'', then the failure of exactness
of the rigidity complex means that Theorem \ref{big theorem}
does not hold.  The degree vanishing and
rigidity (Theorems \ref{degree vanishing}, \ref{scalar automorphisms})
also fail. 

It would be interesting to find a relation between combinatorial
IH and Sabitov's result \cite{Sab} that the volume of a flexible
triangulated $3$-sphere remains constant as it flexes.

\section{New inequalities}
As with any sheaf theory, combinatorial IH sheaves give rise
to a wide variety of homomorphisms, complexes and exact sequences.
We present two applications which produce new inequalities
among the $h$- and $g$-numbers.  The first answers a question
of Stanley from \cite{St87} and the second generalizes Kalai's
monotonicity (Theorem \ref{KM}) specialized at $t=1$.  

\subsection{Decompositions of $h(P,t)$ from a shelling}\label{positive local h}
Let $P$ be a $d$-polytope.
The sum \eqref{h from g} which gives the $g$-polynomial
of $P$ runs over all faces in the boundary
complex $\partial P$, which is a polyhedral subdivision 
of a $(d-1)$-sphere.  In this section, we consider
replacing \eqref{h from g} by 
partial sums over subsets of $\partial P$.  More
specifically, let $I$ be a subcomplex of $\partial P$,
and let $J$ be a subcomplex of $I$.  Then we define
\begin{equation}\label{rel h from g}
h(I,J,t) = \sum_{F\in I\setminus J}g(F,t)(t-1)^{d - 1 - \dim F}.
\end{equation}
This definition is considered in \cite{St87} under the more
general hypothesis that $I$ is a locally Eulerian poset and
$J$ is an order ideal in $I$.

If we have a filtration 
$\emptyset = I_0 \subset I_1 \subset \dots \subset I_r = \partial P$
by subcomplexes, we therefore get a decomposition
\begin{equation}\label{h is sum of local}
h(P,t) = \sum_{j=1}^r h(I_j, I_{j-1}, t).
\end{equation}
If the polynomials $h(I_j, I_{j-1})$ have nonnegative
coefficients, this will result in inequalities for the
$h$-numbers of $P$.  

The following result answers a question of Stanley
\cite[Section 6]{St87}.  Suppose that we are given a shelling
of $\bdy P$; i.e.\ an ordering $F_1,\dots, F_r$ of the
facets of $P$ such that for every $k = 1, \dots, r-1$
the union $I_k = F_1 \cup \dots \cup F_k$
is topologically a $(d-1)$-dimensional disk.
Boundary complexes of polytopes can always be shelled,
for instance by the line shelling
construction of Bruggesser and Mani \cite{BruMa}.

\begin{thm}\label{Shelling polys nonneg}
For $1 \le j\le r$, 
the coefficients of $h(I_j, I_{j-1},t)$ are nonnegative.
\end{thm}

Note that this polynomial is ``locally defined", in 
the sense that it depends only on $F_j$ and 
its intersection with $\bigcup_{i < j} F_i$.
Thus Theorem \ref{Shelling polys nonneg} implies 
that the $h$-numbers are nonnegative for shellable
polyhedral complexes.

It is elementary to see that the sum of the coefficients
of $h(I_j,I_{j-1},t)$ is nonnegative, since setting
$t=1$ in \eqref{rel h from g} 
gives $h(I_j, I_{j-1},1) = g(F_j,1)$.
One case where it is easy to see that
the individual coefficients are nonnegative
is when $F_j$ is a simplex, 
which gives $h(I_j,I_{j-1},t) = t^k$, 
where $k$ is the number of facets of $F_j$ not in 
$\bigcup_{i<j} F_i$. Another is when
$j = 1$ and $j = r$, which gives
\begin{align*}
h(I_1,\emptyset , t) & = t^d g(F_1, t^{-1}) \\
h(I_r, I_{r-1}, t) & = g(F_r, t).
\end{align*}
The nonnegativity for other $j$ is new, although
Bayer \cite{Bay05} showed nonnegativity for 
certain shellings of a class of nonsimplicial 
polytopes generalizing cyclic polytopes.

\begin{example} Let $P$ be a prism over a $2$-simplex, 
and take any shelling
$F_1, \dots, F_5$ for which $F_1$ and $F_4$ are the simplicial 
facets.  Then
\begin{align*}
h(I_1, \emptyset, t) & = t^3 \\
h(I_2, I_1, t) & = 2t^2 \\
h(I_3, I_2, t) & = t + t^2 \\
h(I_4, I_3, t) & = t \\
h(I_5, I_4, t) & = 1 + t
\end{align*}
and $h(P, t) = t^3 + 2t^2 + (t + t^2) + t + (1+t) = 1 + 3t + 3t^2 + t^3$.
\end{example}
  
The proof of Theorem \ref{Shelling polys nonneg} 
relies on the following beautiful result of Barthel, 
Brasselet, Fieseler, and Kaup,
which characterizes when $\cL(\Delta)$ 
is a free $A$-module.  Let $\Delta$ be
a purely $d$-dimensional fan in a $d$-dimensional vector space, and 
let $\partial\Delta$ denote the subfan generated by the $(d-1)$-dimensional
cones which are contained in exactly one $d$-dimensional cone of $\Delta$.
\begin{thm}[\cite{BBFK}]\label{quasi-convex}
$\cL(\Delta)$ is a free $A$-module if and only
if $\partial \Delta$ is an $\R$-homology manifold.  If this holds,
then 
\[\Hilb(\ol{\cL(\Delta)},t) = \sum_{\sig\in \Delta\setminus \bdy\Delta}
g(\sig,t)(t-1)^{d-\dim\sig},\]
where we put $g(\sig, t) = g(F, t)$ when $\sig = cF$, the cone
over a polytope $F$.
\end{thm}
Following \cite{BBFK}, we call fans satisfying 
these equivalent conditions {\em quasi-convex}.  Examples of
quasi-convex fans include complete fans (where $\bdy\Delta$ is
empty) and full-dimensional cones (where $\bdy\Delta$ is 
homeomorphic to $\R^{d-1}$).

\begin{proof}[Proof of Theorem \ref{Shelling polys nonneg}] 
Let $\Delta = \Delta_P$ be the central fan of the $d$-polytope 
$P$, and let $\sig_1,\dots,\sig_r$ be the ordering on the top-dimensional
cones of $\Delta$ given by taking cones over the $F_i$. 

Let $\Delta^j = [\sig_j]\cup \dots\cup [\sig_r]$ be the fan
generated by $\sig_j,\dots,\sig_r$.  
Then $\partial\Delta^j$ is the cone over $\bdy(F_j\cup\dots\cup F_s) = 
\bdy(I_{j-1})$, which is a $(d-1)$-sphere since
the ordering of the $F_i$ is a shelling.  Thus $\Delta^j$ is
quasi-convex.  By Theorem \ref{quasi-convex}, we have
\[\Hilb(\ol{\cL(\Delta^j)},t) = \sum_{F \notin I_{j-1}} 
g(F,t)(t-1)^{d-1-\dim F}.\]

Now for $1\le j\le r$ consider the restriction 
$\cL(\Delta^j)\to \cL(\Delta^{j+1})$, where we set $\Delta^{r+1} = \emptyset$.
It is a surjective map of free $A$-modules, which implies that its kernel 
$K$ is
also free. It follows that
\begin{align*}
\Hilb(\ol{K},t) & = \Hilb(\ol{\cL(\Delta^j)},t) - 
\Hilb(\ol{\cL(\Delta^{j+1})},t) \\
& = \sum_{F \in I_j \setminus I_{j-1}} g(F,t)(t-1)^{d-1-\dim F} \\
& = h(I_j, I_{j-1}, t).
\end{align*}
Thus the coefficients of $h(I_j, I_{j-1}, t)$ are nonnegative.
\end{proof}

\begin{rmk} This result can also be deduced by applying Proposition 6.7 
of \cite{BBFK} 
to the pair of fans $([\sig_i],[\sig_i]\cap \Delta_{i-1})$.
\end{rmk}

When $\Delta = \Delta_P$ for a rational polytope $P$ in $V$, 
and the shelling of $P$ is a line shelling, 
the local $h$-polynomials can be understood in terms
of a topological construction.  The line shelling 
depends on the choice of a vector $v \in V$; if this
vector is in the lattice, it determines a homomorphism
from $\C^*\to T$ and thus an action 
of $\C^*$ on the toric variety $X_\Delta$.  The condition
that this direction determines a shelling ensures that 
the action has isolated fixed points.  

The resulting flow on $X_\Delta$ induces a partial order on the set
$X_\Delta^T$ of fixed points.  The fixed points
are in bijection with the facets of $P$, and the 
partial order is compatible with the total order from 
the shelling.  The decomposition \eqref{h is sum of local}
then comes from a result of Kirwan \cite{Kir}
which gives a decomposition of the equivariant 
intersection cohomology of a $T$-variety $X$ 
into a sum of terms coming from each component of the
locus of points fixed by a one-dimensional subtorus of $T$.

\subsection{A generalization of Kalai's monotonicity at $t=1$} \label{generalized monotonicity}
Our next application gives a lower bound for $g(P, 1) = \sum g_i(P)$
 which generalizes Theorem \ref{KM}
specialized to $t = 1$.  For the following discussion it will be
convenient to work entirely with cones instead of polytopes;
if $\sig$ is the cone $cP$ over a compact polytope $P$,
we put $g(\sig,t) = g(P,t)$.  

Fix a full-dimensional cone $\sig$ in a vector space $V$, and let 
$\Delta = [\sig]$. Fix a vector $v \in V$, and let $\Delta_0$ denote
the set of faces $\tau \in \Delta$ for which $v \in \Span(\tau)$.

\begin{thm} \label{generalized Kalai} Let $\Min \Delta_0$ 
be the set of minimal faces in $\Delta_0$.  Then
\begin{equation} \label{ineq}
g(\sig,1) \ge \sum_{\tau \in \Min \Delta_0} g(\tau,1)g(\sig/\tau,1).
\end{equation}
\end{thm}

Unlike the other results described so far, this is not a purely
combinatorial statement about face posets of cones and polytopes,
since in general the possible sets $\Min \Delta_0$ will 
depend on the angles of the faces with respect to $v$, and not 
just their inclusion relations. There is one case where 
$\Delta_0$ can be determined combinatorially, however: 
when $v$ lies in the relative interior 
of a face $\tau$ of $\sig$, then $\Delta_0 = \st(\tau)$. 
Therefore
$\Min\Delta_0 = \{\tau\}$, so Theorem \ref{generalized Kalai}
becomes \[g(\sig,1) \ge g(\tau,1)g(\sig/\tau,1).\]
If $\sig = cP$ and $\tau = cF$ for a polytope $P$ and a face 
$F$ of $P$, then $\sig/\tau = c(P/F)$, and we recover Theorem
\ref{KM} at $t = 1$.

\begin{example}
Let $\sig$ be the cone over a square, and choose $v$ 
so that $\Delta_0$ contains two opposite $2$-faces $\tau_1$, $\tau_2$
in addition to $\sig$ itself.  Then $\Min \Delta_0 = \{\tau_1, \tau_2\}$,
and Theorem \ref{generalized Kalai} says that
\[g(\sig, 1) \ge \sum_{i=1,2} g(\tau_i,1)g(\sig/\tau_i,1) = 
1\cdot 1 + 1 \cdot 1 = 2.\]
In this case equality holds, since $g(\sig,t) = 1 + t$.  This example
shows that evaluating at $t = 1$ is necessary. 
\end{example}

To prove Theorem \ref{generalized Kalai}, we restrict the combinatorial
IH sheaf $\cL_\Delta$ from $\Delta$ to $\Delta_0$ in a 
way which uses the direction of the vector $v$.  It is a combinatorial
counterpart of the ``hyperbolic localization" of \cite{Br03}.
In order to describe it, decompose the fan $\Delta = [\sig]$ according
to its intersection with lines parallel to $v$.

\begin{defn}
Call a face $\tau \in \Delta$ a {\em back} face if  
for every  point $x \in \tau^\circ$,
the set 
\[\{t\in \R \mid x + tv\in \sig\}\]
contains some interval $[0,\epsilon_x), \epsilon_x>0$.
Let $\Delta_{\le 0}$ be the set of all back faces.
Using the same definitions but with $v$ replaced by 
$-v$, define a subset $\Delta_{\ge 0}\subset \Delta$ of {\em front}
faces.
\end{defn}

It is easy to see that
$\Delta_0 = \Delta_{\le 0} \cap \Delta_{\ge 0}$.
We also put
$\Delta_+ = \Delta_{\ge 0} \setminus \Delta_0$, and 
$\Delta_- = \Delta_{\le 0} \setminus \Delta_0$.

If $\tau\in \Delta_{\le 0}$ is a back face, then there is a 
unique face $\rho\prec \sig$ so that
$x + tv\in \rho^\circ$ 
for any $x\in \tau^\circ$ and 
any $t\in (0,\epsilon_x)$; we
denote this face by $\tau_{+v}$.  Note that
if $\tau \in \Delta_0$, then $\tau_{+v} = \tau$.

\begin{prop} \label{front and back}
These sets satisfy the following:
\begin{enumerate}
\item $\Delta_{\ge 0}$, $\Delta_{\le 0}$, and $\Delta_0$ are closed
in the fan topology on $\Delta$.
\item $\tau\mapsto \tau_{+v}$ is a surjective function 
$\Delta_{\le 0}\to \Delta_0$;
$\tau_{+v}$ is the unique smallest face in $\Delta_0$ with
$\tau\prec \tau_{+v}$.
\item $[\Delta_+] = \Delta \setminus \Delta_{\le 0}$
\item The projection $p\colon V \to V/\R v$ gives a 
conewise linear isomorphism between 
$[\Delta_+]$ and a fan $p_*[\Delta_+]$
with support $p(|\Delta|) = p(\sig)$.
\end{enumerate}
\end{prop}

The last two statements can be understood as follows: the subfan 
$[\Delta_+]$ consists of all faces of $\sig$ illuminated by light
shining from infinity in parallel beams with direction $-v$. 

Given an arbitrary subset $\Sig \subset \Delta$, define 
$\Sig_{\ge 0} = \Delta_{\ge 0} \cap \Sig$ and
$\Sig_+ = \Delta_+ \cap \Sig$. 

\begin{defn}
Given a $\cA_\Delta$-sheaf $\cF$, its $v$-localization is
the sheaf $\cF^v$ defined by  
\[\cF^v(\Sig)= \cF([\Sig_{\ge 0}],[\Sig_+])\]
for any subfan $\Sig$ of $\Delta$.
\end{defn}  
The support of $\cF^v$ is contained in $\Delta_0$, since
if $\Sig\cap \Delta_0 = \emptyset$, then $\Sig_{\ge 0} = \Sig_+$, and so
$\cF^v(\Sig) = 0$.

\begin{thm} \label{purity} If $\cF$ is a pure $\cA_\Delta$-sheaf, then 
$\cF^v$ is also pure.
\end{thm}

\begin{prop} \label{local Euler char} 
For every $\tau\in \Delta_0$, we have
$\dim_\R \ol{\cF^v(\tau)}= \dim_\R \ol{\cF(\tau)}$.
\end{prop}
Note that the Hilbert polynomials of these modules are
not in general equal; it is essential here to ignore the grading and
evaluate at $t=1$.

We now deduce Theorem \ref{generalized Kalai} from these results.  
Let $\cL = {}_{\zcone}\cL$.
By Theorem \ref{purity} $\cL^v$ is a direct sum of
shifted copies of sheaves ${}_\tau\cL$, $\tau\in \Delta$.  Since
the support of $\cL^v$ lies in $\Delta_0$, we have
\begin{equation}
\label{xyz} 
\cL^v \cong \bigoplus_{\tau \in \Delta_0} {}_\tau\cL \otimes N_\tau,
\end{equation}
where the $N_\tau$ are graded $\R$-vector spaces.

Proposition \ref{local Euler char} implies that
$\dim_\R N_\tau = \dim_\R \ol{\cL(\tau)} = g(\tau, 1)$ if
$\tau\in \Min \Delta_0$.  By taking just those terms in \eqref{xyz} 
with $\tau\in \Min \Delta_0$, we get
\[g(\sig, 1) = \dim_\R \cL^v(\tau) \ge \sum_{\tau\in \Min\Delta_0}
\dim_\R (\ol{{}_\tau\cL(\sig)})\dim_\R (N_\tau) = 
\sum_{\tau\in \Min\Delta_0} g(\tau, 1)g(\sig/\tau, 1),\]
which is Theorem \ref{generalized Kalai}.

Note that when $v$ lies in the relative interior of a cone
$\tau\prec\sig$, we have $\Delta_{\ge 0} = \st(\tau)$, 
$\Delta_+ = \emptyset$, and so
 $\cL^v = \cL|_{\st(\tau)}$.  The same argument
then works without appealing to Proposition \ref{local Euler char},
and in fact gives a graded statement without evaluating 
at $t=1$.  There is only one term in \eqref{ineq}, since $\Min\Delta_0 = \{\tau\}$.
This is precisely the argument 
in \cite{BBFK,BrLu} used to
prove Theorem \ref{KM}.

\subsection{Proofs} We give the proofs of 
Theorem \ref{purity} and Proposition \ref{local Euler char}.

First we show that if $\cF$ is a flabby sheaf, then $\cF^v$ 
is again a flabby sheaf.  We use the following characterization
of flabby sheaves.  Note that flabbiness of an $\cA_\Delta$-module
only depends on its structure as a sheaf of graded vector spaces, and
not the module structure.  For a cone $\tau\in \Delta$, let 
$\Ru_{\st(\tau)}$ be the constant sheaf on $\st(\tau)$
with stalk $\R$, extended by zero to all of $\Delta$.
Then $\Ru_{\st(\tau)}(\Sig) = \R$ if
$\tau\in \Sig$ and is $0$ if $\tau\notin\Sig$, and the restriction
maps are the identity whenever possible.   The following lemma is
essentially Lemma 3.6 from \cite{BBFK}. 
\begin{lemma} 
\label{injective=flabby}
A sheaf of graded vector spaces on $\Delta$ 
is flabby if and only if it is
a direct sum of sheaves $\Ru_{\st(\tau)}$, $\tau\in\Delta$, 
with shifts.
\end{lemma}

Using this, the following lemma shows that if $\cF$ is flabby then
$\cF^v$ is flabby.

\begin{lemma} The $v$-localization of $\cF = \Ru_{\st(\tau)}$ is given by 
\[\cF^v \cong \left\{\begin{array}{ll} \Ru_{\st(\tau_{+v})}, & 
\tau\in \Delta_{\le 0} \\
0 & \text{otherwise.} 
\end{array}\right.\]
\end{lemma}
\begin{proof} If $\Sig$ is a subfan of $\Delta$, then
\[\cF^v(\Sig) = \left\{\begin{array}{ll} \R, & 
\tau\in [\Sig_{\ge 0}]\setminus [\Sig_+] \\
0 & \text{otherwise.} 
\end{array}\right.\]
Suppose first that $\tau\in [\Sig_{\ge 0}]\setminus [\Sig_+]$ .  Then there exists
$\rho\in \Delta_{\ge 0} \setminus \Delta_+ = \Delta_0$ with $\tau\prec\rho$.
If $\tau$ were in $[\Delta_+]$
then there would be a cone $\nu \in \Delta_+$ with $\tau\prec \nu\prec \rho$.
This is because the projection $p(\rho)$ is a union of cones containing 
$\tau$ in the
projected fan $p_*[\Delta_+]$, so one can take $\nu$ to 
be a cone which projects to a
maximal-dimensional cone contained in $p(\rho)$ and containing $p(\tau)$.
But this means that $\nu\in[\rho] \subset \Sig$, so $\tau\in [\Sig_+]$, contrary to
assumption.  So $\tau\in \Delta_0\setminus [\Delta_+] = \Delta_{\le 0}$.  
Then part (2) of Proposition \ref{front and back} implies that 
$\tau_{+v} \prec \rho$, so $\Sig \cap \st(\tau_{+v}) \ne\emptyset$.

Conversely, if $\tau \in \Delta_{\le 0}$ and 
$\Sig \cap \st(\tau_{+v}) \ne\emptyset$, then in 
particular $\tau_{+v} \in \Sig\cap\Delta_0 = \Sig_0$,
so $\tau \in [\Sig_{\ge 0}] \setminus [\Sig_+]$.
\end{proof}

Next we show that $\cF^v$ is locally free if $\cF$ is pure.  
Take a cone $\tau\in \Delta_0$.  Then $[\tau]_{\ge 0} = [\tau]\cap \Delta_{\ge 0}$
is the set of front faces of $\tau$, considered as a cone in
$\Span(\tau)$, with respect to the vector $v$, 
and similarly for the back faces.
So to show that $\cF^v$ is locally free, it is enough to show that the
stalk at $\sig$
\[\cF^v([\sig]) = \cF([\Delta_{\ge 0}], [\Delta_+])\]
is a free $A$-module.  This will follow from the following result
and Proposition \ref{front and back}(4).

\begin{thm}
Let $\sig$ be a full-dimensional cone in $V$, and 
let $\cF$ be a pure $\cA_{[\sig]}$-module.
Suppose that $\Sig$ is a purely $(d-1)$-dimensional subfan of 
$\bdy\sig$ which is conewise linearly 
isomorphic to a quasi-convex fan via the projection
$p\colon V \to V/\R v$ for some $v\in V$.  Then
\item $\cF([\sig],\Sig)$ is a free $A$-module, of the
same total rank as $\cF(\sig)$.
\end{thm}
\begin{proof}
An $A$-module $M$ is free if and only if
$\Tor^{A}_i(M,\R) = 0$ for $i > 0$ (see \cite[\S0.B]{BBFK}).  
Using the long exact Tor sequence of
\[0 \to \cF([\sig],\Sig)  \to \cF([\sig])\to \cF(\Sig)\to 0\]
and the freeness of $\cF([\sig])$, we are reduced to showing 
that $\Tor^A_i(\cF(\Sig),\R) = 0$ for $i > 1$
and that $\dim \Tor_0^A(\cF(\Sig)) = \dim \Tor_1^A(\cF(\Sig))$.
  
To see this, let $\hat{\Sig}$ denote the fan in $\hat V = V/\R v$
obtained by projecting $\Sig$ by $p$.  Also 
let $\hat{A}$ denote the ring of polynomial functions on 
$V/\R v$; the pullback $p^*\colon \hat A\to A$ makes it
into a subring of $A$.  

Since $p$ is a conewise linear isomorphism, pushing forward
the restriction of $\cF$ to ${\Sig}$ gives 
a pure sheaf on $\hat\Sig$ whose global sections are  
just $\cF(\Sig)$, with the $\hat A$-module structure induced
by $p^*$.  This means that $\cF(\Sig)$ is a free $\hat A$-module,
by Theorem \ref{quasi-convex}.

Let $M = \cF(\Sig)$.  Then the $A$-module structure is determined
by the $\hat A$-module structure together with the homomorphism
$M\to M$ of multiplication by $y$, where $y$ is any
element in $A_1$ not in $\hat A_1$.  

There is a short exact sequence
\[0 \to A \otimes_{\hat A} M \to A \otimes_{\hat A} M \to M\to 0\] 
of $A$-modules, where 
the first map is $a \otimes m \mapsto a \otimes ym - ya \otimes m$  and 
the second map is $a \otimes m \mapsto am$.  The 
first two terms of this sequence give a free resolution of $M$ as an $A$-module, and tensoring
with $\R$ gives a finite-dimensional 
two-step complex $M \otimes_{\hat A} \R \to M \otimes_{\hat A} \R$
whose homology in degree $i$ is $\Tor^A_i(M, \R)$.  The theorem follows.
\end{proof}

\section{Polytope duality and Stanley's convolution identity}\label{Koszul}
As a final application, we explain how an identity of Stanley
relating the $g$-numbers of a polytope $P$ and its polar 
$P^*$ can be understood in the combinatorial IH language.
We also discuss how this is related to the toric Koszul 
duality constructed in \cite{BraLu}.

For a polytope $P$ which spans the linear space $V$, the polar
polytope is given by
\[P^* = \{w \in V^* \mid \la v , w \ra \ge -1\; \text{for all}\; v\in P\},\]
where $P$ is translated so that $0$ lies in its interior;
it is well-defined up to projective equivalence.  There is an 
order-reversing inclusion between the face lattices of 
$P$ and $P^*$.  

In \cite{St92}, Stanley showed that for 
any polytope $P \ne \emptyset$ we have
\begin{equation}\label{reciprocity} 
\sum_{\emptyset \le F \le P} (-1)^{\dim F}g(F^*,t)g(P/F,t) = 0.
\end{equation} 
In particular, if $\dim P = 2k$ is even, then the degree $k$ piece of 
this identity gives $g_k(P) = g_k(P^*)$.  In fact 
\eqref{reciprocity} holds for general Eulerian posets, but our
interpretation will only be valid for polytopes.

We will show that Stanley's formula can be
``lifted'' to a statement in linear algebra, by exhibiting 
a long exact sequence 
of graded vector spaces whose graded Euler characteristic is the 
alternating sum in \eqref{reciprocity}. It is obtained by 
taking stalk cohomology of a certain complex  
of $\cA_\Delta$-modules.  A more detailed discussion of the homological
algebra of $\cA_\Delta$-modules appears in \cite[Section 6]{BraLu},
but for our purposes it is enough to remark that a sequence
$0 \to \cE \to \cF \to \cG\to 0$
of $\cA_\Delta$-modules is short exact if and only if the induced
sequence on stalks
$0 \to \cE(\sig) \to \cF(\sig) \to \cG(\sig)\to 0$
is a short exact sequence of $A_\sig$-modules for every $\sig\in\Delta$. 
 
The complex we want is described by the following theorem.
Fix a fan $\Delta$ and let $r = \dim \Delta$ be the dimension of 
its largest cone. Let $\cM$ be the $\cA_\Delta$-module which is
the extension by zero of the rank one constant sheaf on $\{o\}$;
in other words, we have $\cM(o) = \R$ and $\cM(\tau) = 0$ for $\tau \ne o$.
\begin{thm}\label{Verma complex} 
There exists a resolution 
\begin{equation} \label{Verma resolution}
0 \to \cM \to \cF^0 \to \cF^1 \to \dots \to \cF^r \to 0
\end{equation}
of $\cM$ by pure sheaves so that $\cF^0 = {}_o\cL_\Delta$ and
for $i > 0$ $\cF^i$ is a direct sum of sheaves of 
the form ${}_\tau\cL_\Delta[-\frac{\dim \tau - i}{2}]$,
with $\dim \tau - i \in 2\Z_{\ge 0}$; it is unique up 
to a unique automorphism, after fixing a basis of 
$\cM(0) = \R$.
\end{thm}

We postpone the proof until the end of this section.

We are interested in the multiplicity $m_k(\tau)$ of the 
sheaf ${}_\tau\cL_\Delta[-k]$ in 
$\cF^{\dim\tau - 2k}$.  Since the complex provided by 
Theorem \ref{Verma complex} clearly behaves well under 
restriction of fans, this multiplicity depends only on $\tau$
and not on the ambient fan $\Delta$.
It can be computed recursively as follows.
If $\tau\ne o$, the stalk complex $\cF^\udot(\tau)$ is
an exact sequence of free $A_\tau$-modules, so the
reduced sequence
\begin{equation}\label{reduced sequence}
0\to \ol{\cF^0(\tau)} \to \dots \to \ol{\cF^r(\tau)} \to 0
\end{equation}
is exact.  Taking the Euler characteristic of the graded 
piece in degree $k$, we can solve for  
$m_k(\tau)$, given the multiplicities $m_j(\rho)$ for all
$j \le k$ and all proper faces $\rho$ of $\tau$.  Using
\eqref{reciprocity}, we get the following result.
\begin{prop} \label{Verma multiplicities}
If $\tau$ is a cone over a polytope $P$, then
$m_k(\tau) = g_{k}(P^*)$
for all $k\ge 0$.  
\end{prop}
Thus the graded Euler characteristic of the exact sequence
\eqref{reduced sequence} gives Stanley's formula.
The dimensions of the individual entries can also 
be organized into a two-variable polynomial $B(P;u,v)$
considered in \cite{BatBor,BorMav} as part of their 
study of stringy Hodge numbers for hypersurfaces in
toric varieties.   

As usual, an exact sequence gives rise to inequalities
by truncation.
Taking the first $s+1$ terms of the degree $k$ piece of 
\eqref{reduced sequence}, we see that
\begin{equation} \label{truncated inequality}
\sum (-1)^{\dim F - s + 1} g_i(F^*)g_j(P/F) \ge 0,
\end{equation}
where the sum is over all $i$ and $j$ with $i + j = k$ and 
all faces $\emptyset \le F \le P$ with $\dim F \le s + 2i-1$.
Taking $s = 0$ just recovers the inequality $g_k(P)\ge 0$, since
the only nonzero summand in \eqref{truncated inequality} is 
$F = \emptyset$.  For $k = s = 1$, the inequality can be reduced to
to $g_2 + f_1 - f_0 \ge (d+1)(d-2)/2$, where $d = \dim P$.  This
follows from the nonnegativity of $g_2$ and 
the elementary facts that $f_0 \ge d+1$ and every vertex is 
contained in at least $d$ edges.  For $d \le 5$ all the
other cases of \eqref{truncated inequality} follow from these 
two cases.  In general we do not know if any of these inequalities 
are new in dimensions $d \ge 6$.
 
Since the complex $\cF^\udot$ is unique up to a unique 
isomorphism, we can define a canonical 
multiplicity space of dimension $m_k(\tau)$; let
$M_k(\tau)$ be the image of the degree $k$
part of the map $\ol{\cF^{i}(\tau,\bdy\tau)}\to 
\ol{\cF^{i}(\tau)}$, where $i = \dim(\tau) - 2k$. 
Then we have a canonical isomorphism
\[\cF^i \cong \bigoplus_{\dim\tau = 2k + i} M_k(\tau)\otimes_\R {}_\tau\cL[-k].\]
Proposition \ref{Verma multiplicities} then
lifts to the following functorial statement.  Recall that 
the dual cone to a full-dimensional cone $\tau$ in $V$ is 
\[\tau^\vee = \{y \in V^* \mid \la x, y\ra \ge 0\; \text{for all}\; x\in \tau\}.\]
If $\tau$ is the cone over a polytope $P$, then $\tau^\vee$ is the cone over $P^*$.
\begin{thm}\label{canonical}
There is a canonical isomorphism
\[M_k(\tau)\cong \ori(\tau)\otimes_\R\ol{\cL(\tau^\vee)}_k^{\,*},\]
where $\ori(\tau)\cong \R$ is the space of
orientations of $\tau$.  
\end{thm}

There is an appealing special case of this result when
$\dim \tau = 2k+1$.  The degree $k$ part of \eqref{reduced sequence}
gives an isomorphism $\ol{\cF^0(\tau)}_k \cong M_k(\tau)$.  
Thus there is a canonical dual pairing
\[\ol{\cL(\tau)}_k \otimes \ol{\cL(\tau^\vee)}_k \to \ori(\tau) \cong
\ori(\tau^\vee).\]
It is a pleasant exercise in linear algebra to construct such a 
pairing for $k = 1$ using the descriptions of these vector spaces
in Section \ref{degree zero and one}.  
Giving an explicit description
of such a pairing for $k=2$ in terms of the discussion in 
Section \ref{degree two} seems to be much more difficult.

\subsection{} We sketch the proof of Theorem \ref{canonical} 
using the toric Koszul duality constructed in \cite{BraLu}.  
In that paper, it is shown that the derived category $D^b(\cA_\Delta)$
of $\cA_\Delta$-modules 
is a ``mixed" version of the topological equivariant derived category 
$D^b_T(X_\Delta)$ of the toric variety $X_\Delta$.  Roughly, this means
that $D^b_T(X_\Delta)$ can be expressed as the derived category of modules
over an associative ring $R$, and $D^b(\cA_\Delta)$ is the 
derived category of graded modules over a graded version of $R$; these
two categories are then related by a functor which forgets the grading.
Thus $D^b(\cA_\Delta)$ has two independent shift functors
coming from shifting the grading in the complex and shifting 
the algebraic grading of the modules (note that the grading shift on
graded $R$-modules is not the same as the grading shift on 
$\cA_\Delta$-modules).

Every object in $D^b(\cA_\Delta)$ is isomorphic to a
complex of pure sheaves; pure sheaves and resolutions by them
play a role similar to the role of injectives in other derived
categories of sheaves.  For instance, the resolution given by 
Theorem \ref{Verma complex} gives an isomorphism 
$\cM \stackrel{\sim}{\to}\cF^\udot$ in $D^b(\cA_\Delta)$.
The degree restrictions on the summands of the $\cF^i$
says that this object $\cM$ lies in the abelian subcategory  
$P(\cA_\Delta) \subset D^b(\cA)$ of perverse objects,
which corresponds to the abelian category of graded $R$-modules.
The corresponding topological objects are
equivariant perverse sheaves on $X_\Delta$.  For instance, $\cM$ 
represents the ``Verma'' or ``standard'' sheaf on 
$X_\Delta$ which is the extension by $0$ of the rank one 
constant local system on the open orbit; this is perverse since
the inclusion of the orbit is an affine map. 

The simple objects in $P(\cA_\Delta)$ are the sheaves
${}_\tau\cL$, taken with certain shifts which we ignore.
The combinatorial IH sheaves in the complex $\cF^\udot$ 
are the simple constituents of the perverse object $\cM$,
and the ``dumb'' filtration obtained by truncating 
the complex gives the weight filtration on $\cM$, which has
semi-simple subquotients. 

There are enough projectives in $P(\cA_\Delta)$, so we can take a 
projective cover ${}_\tau\cP^\udot \to {}_\tau\cL$.
The total multiplicity space $M_\udot(\tau) = \bigoplus_i M_i(\tau)$
is canonically isomorphic to the total Hom-space
$\underline{\Hom}({}_\tau\cP^\udot, \cM)$
which takes all homomorphisms with all possible grading
shifts.

Now suppose that $\Delta = [\tau]$ for a full-dimensional rational cone 
$\tau$, and let $\Delta^\vee = [\tau^\vee]$.  
Then the main result of \cite{BraLu} constructs an equivalence of triangulated 
categories \begin{equation*}
K\colon D^b(\cA_\Delta) \to D^b(LC_\cF(X_{\Delta^\vee}))
\end{equation*}
called the Koszul duality functor.  The category on the
right hand side is the derived category of
sheaves on $X_{\Delta^\vee}$ 
which are constructible with respect to the orbit 
stratification on $X_{\Delta^\vee}$ and are endowed with an 
extra ``mixed'' structure.

Applying the Koszul functor $K$ sends ${}_\tau\cP$ to the
intersection cohomology complex $IC^\udot(X_{\Delta^\vee})$,
appropriately incarnated as an object in 
$D^b(LC_\cF(X_{\Delta^\vee}))$.  The image of $\cM$ under $K$
is isomorphic to the 
point sheaf $i_*\R_p$ where $\{p\}$ is the the unique torus
fixed point of $X_{\Delta^\vee}$, and $i\colon \{p\} \to
X_{\Delta^\vee}$ is the inclusion.
 
Since $K$ is an equivalence of categories, it induces an 
isomorphism
\[\underline{\Hom}({}_\tau\cP^\udot, \cM)
 \stackrel{K}{\lra} 
\underline{\Hom}(IC^\udot(X_{\Delta^\vee}), i_*\R_p).\] 
The left side is identified with $M_\udot(\tau)$, as we explained above,
while the right side is the dual to the stalk intersection 
cohomology of $X_{\Delta^\vee}$ by a standard adjunction.
This essentially proves Theorem \ref{canonical}.  The 
fact that the graded pieces on each side correspond
follows from how the functor $K$ behaves under the
shifts and twists from the mixed structure.  The
appearance of the orientation group in the theorem comes
from pinning down $K(\cM)$ more precisely; it turns out to 
be {\em canonically} isomorphic to $i_*(\ori(\tau)_p)$.  
This can be seen by noticing that $M_0(\tau) \cong \ori(\tau)$,
which follows from the proof of Lemma \ref{a lemma} below.

\subsection{Proof of Theorem \ref{Verma complex}}

We proceed by induction on the number of cones in the fan $\Delta$.
If there is only one cone, then $\Delta = \{o\}$ and the result is 
obvious.  If $\Delta$ has more than one cone, let 
$\sig\in \Delta$ be a maximal cone, and
let $\wt\Delta = \Delta \setminus \{\sig\}$.  Assume 
inductively that we have constructed the required resolution 
$0 \to \cM \to \wt\cF^0 \to \wt\cF^1 \to \dots \to \wt\cF^s \to 0$,
$s = \dim \wt\Delta$ of sheaves
on $\wt\Delta$.  We will extend it to a resolution on $\Delta$.

For each $i$, let $\cE^i$ be the sheaf on $\Delta$ which is
the minimal pure extension of $\wt\cF^i$.  In
other words, $\cE^i|_{\wt\Delta} = \wt\cF^i$, 
and the restriction map makes $\cE^i(\sig)$ into the 
minimal free cover of $\cE^i(\bdy\sig)$.  
Using the degree bounds from Theorem \ref{degree vanishing}, it is 
easy to see that
the maps $\wt\cF^i \to \wt\cF^{i+1}$ extend uniquely to maps
$\cE^i\to \cE^{i+1}$ (see \cite[Theorem 6.6.2]{BraLu}).  
The resulting sequence of sheaves and maps
\[0 \to \cM \to \cE^0 \stackrel{d^1}\lra \cE^1 \stackrel{d^2}\lra \,\cdots\, \stackrel{d^s}\lra \cE^s\]
is almost the resolution we want, but the composition of successive maps
may not be zero, and $d^s$ may not be surjective.
We will add sheaves which are supported
only on the cone $\sig$
in order to recover the chain complex property.

Let $n = \dim\tau$, so $r= \dim\Delta = \max(n,s)$. 
For any $i\ne n$, let $\psi^i\colon \cE^{i-1} \to \cE^{i+1}$ be the 
composition
$d^{i+1}\circ d^i$, and let $\cG^i = \im \psi^i$.  
Since $\cE^\udot|_{\bdy\sig}$
is a complex by assumption, $\cG^i$ has nonzero stalk only on $\sig$. 
Using the degree restrictions on the simple constituents of the 
$\cE^j$ together with Theorem \ref{degree vanishing}, we see that
$\cE^{i-1}(\sig)$ is a free $A_\sig$-module generated in degrees 
$\le \delta(i)$, where we put 
$\delta(j) = (n-j)/2$.
Further, $\cG^i(\sig)$ is contained in $\cE^{i+1}(\sig,\bdy\sig)$, which is
a free $A_\sig$-module generated in degrees $\ge \delta(i)$.  It follows that 
$\cG^i(\sig)$ is a free $A_\sig$-module generated in degree exactly 
$\delta(i)$. 
In particular, $\cG^i = 0$ if $i$ and $n$ have different parity or if 
$i\ge n$. 

Let $\cG^n$ be the cokernel of $d^{n-1}$ on $[\sig]$, extended by zero
to the other cones; it is 
supported only on $\sig$.
Since $\cE^{n-2}(\sig)$ and $\cE^{n-1}(\sig)$ are free $A_\sig$-modules
generated in degree $0$, so is $\cG^n(\sig)$.

Define $\cF^i = \cE^i \oplus \cG^i$. We make these sheaves into 
a complex by defining the boundary maps to be
\[\dots \to \cE^{i-1} \xrightarrow{\genfrac{[}{]}{0pt}{}{d^{i}}{-\psi^i}
} \cE^i \oplus \cG^i \xrightarrow{[d^{i+1}\; \iota]} 
\cE^{i+1} \to \cdots\]
for any $i \ne n$ with $i-n$ even, where $\iota\colon \cG^i \to \cE^{i+1}$ 
is the inclusion,  and defining the boundary maps at the $n$th
position to be
\[\dots \to \cE^{n-1} \xrightarrow{\genfrac{[}{]}{0pt}{}{d^{n}}{p}
} \cE^n \oplus \cG^n \xrightarrow{[d^{n+1}\;0]} 
\cE^{n+1} \to \cdots,\]
where $p\colon \cE^{n-1} \to \coker d^{n-1}$ is the natural projection,
extended by zero.

It is easy to check that this is a complex.  To see that it is
a resolution of $\cM$, we need to show that $\cF^\udot(\sig)$
is an exact sequence, since $\cM(\sig) = 0$.  Since
 $\cF^\udot(\sig)$ is a complex of free $A_\sig$-modules, it
will be exact if and only if the reduced complex
$\ol{\cF^\oudot(\sig)}$ is an exact sequence of graded vector
spaces.  The degree $k$ part $\ol{\cF^i(\sig)}_k$
vanishes if $i > n-2k$, and if $i = n - 2k$ (so $\delta(i) = k$),
the boundary
$\ol{\cF^{i-1}(\sig)}_k \to \ol{\cF^{i}(\sig)}_k = \cG^i(\sig_k)$ 
is surjective by
construction.  So we only need to check exactness at $i < n - 2k$.

\begin{lemma} \label{a lemma}
$H^i(\cF^\udot(\bdy\sig)) = 0$ if $i \ne n - 1$,
and $H^{n-1}(\cF^\udot(\bdy\sig))$ is isomorphic to 
$\R$, placed in degree zero.  
\end{lemma}

Assuming this for the moment, 
consider the short exact sequence
\[ 0\to \cF^\udot(\sig,\bdy\sig)\to \cF^\udot(\sig) \to 
\cF^\udot(\bdy\sig)\to 0
\]
of chain complexes.
It induces a long exact cohomology sequence
\[\dots \to H^i(\ol{\cF^\oudot(\sig,\bdy\sig)}) \to
H^i(\ol{\cF^\oudot(\sig)}) \to 
H^i(\cF^\udot(\bdy\sig) \stackrel{L}\otimes_{A_\sig} \R)\to \cdots,\] 
where the third term is the cohomology of a derived functor, 
since $\cF^\udot(\bdy\sig)$ is not a complex of free modules.
The first term vanishes in degrees $< \delta(i)$, the
second vanishes in degrees $> \delta(i)$, and the
last term is nonzero only in degree $n-1-i$, since by the
lemma $\cF^\udot(\bdy\sig)$ has a free resolution by a
Koszul complex.

Combined with our previous observations, this
implies the vanishing of $H^i(\ol{\cF^{\oudot}(\sig)})_k$ 
except possibly when $k = 0$ and $i=n-1$.  
But the complexes $\cF^\udot(\sig)^{}_0$
and $\cF^\udot(\bdy\sig)^{}_0$ are the same except 
that $\cF^k(\sig)_0 = \R$ while $\cF^k(\bdy\sig)^{}_0 = 0$.
Lemma \ref{a lemma} thus implies that the Euler characteristic
of $\cF^\udot(\sig)^{}_0$ is zero, so the last cohomology
group also vanishes, completing the construction of 
our resolution.

The uniqueness of $\cF^\udot$ is easily proved by 
following the same 
induction.  Any automorphism of the complex $\wt{\cF^\udot}$ 
extends uniquely to an automorphism of $\cE^\udot$, since the 
sheaves ${}_\tau\cL$ themselves have only scalar automorphisms.  
This extends uniquely to an automorphism of
$\cF^\udot$, since 
the maps $\psi^i\colon \cE^{i-1} \to \cG^i$ and
$p \colon \cE^{n-1} \to \cG^n$ are surjective.  

\begin{proof}[Proof of Lemma \ref{a lemma}]
The global sections of an exact sequence of flabby sheaves
is exact, so the positive degree parts of the sequence 
$\cF^\udot(\bdy\sig)$ are exact sequences of vector spaces.

In degree zero, it is easy to see by induction that
$\cF^i_0 = \bigoplus_{\dim \tau = i} \Ru_{\st(\tau)}$,
and that the component of the 
boundary map $\Ru_{\st(\tau)} \to\Ru_{\st(\rho)}$ for
$\dim \rho = \dim \tau + 1$ is non-zero if and only if 
$\tau$ is a face of $\rho$.  Taking global sections, the complex
$\cF^\udot(\bdy\sig)_0$ is isomorphic to the augmented cellular 
chain complex of the regular cell complex obtained by intersecting the 
fan $\bdy\sig$ with a sphere centered at the origin.
This cell complex is a $(n-1)$-sphere.  
\end{proof}  

\section*{Appendix: two results of Kalai}
With his kind permission, we present Kalai's previously 
unpublished proofs of two consequences
of his monotonicity conjecture (Theorem \ref{KM}).  
The first relates the $g$-numbers of a polytope $P$ and its polar $P^*$.
\begin{thm}[Kalai]
If $g_k(P) = 0$, then $g_k(P^*) = 0$.
\end{thm}
\begin{proof} Use induction on $k$ and $d$.  
The case $k = 0$ is trivial, since $g_0(P) = 1$ for
all polytopes $P$.  Similarly, 
if $d < 2r$, then $g_r(P) = g_r(P^*) = 0$,
and there is nothing to prove. 

Otherwise, assume the result holds when $k < r$ or $k=r$ and $\dim P < d$.  
The degree $r$ term of Stanley's formula \eqref{reciprocity} gives 
\begin{equation}\label{Reciprocity degree k}
\sum_{\emptyset \le F \le P}\sum_{i+j = r} (-1)^{\dim F}g^{}_i(F^*)g_j(P/F) = 0.
\end{equation}
If $g_k(P) = 0$, then Theorem \ref{KM} implies that $g^{}_i(F)g_j(P/F) = 0$ for all
$i + j = k$ and all faces $F$ of $P$.  By the inductive hypothesis this 
means that all the terms of \eqref{Reciprocity degree k} vanish except for $g_k(P^*)g_0(P/P) = g_k(P^*)$.
\end{proof}

The next result is Theorem \ref{Kalai's theorem} from Section 1, which 
says that $g_k(P) = 0$ implies $g_{k+1}(P) = 0$.
This follows from Theorem \ref{KM} and the following identity.
\begin{prop}[Kalai] \label{Kalai's prop}
For any $d$-polytope $P$ and $0\le k \le d/2 - 1$,
\begin{equation*} 
(k+1)g_{k+1}(P) + (d-k+1)g_{k}(P) = 
       \sum_{i=0}^k (i+1)\!\!\!\mathop{\sum_{F \le P}}_{\dim F = 2i}\!\! g_i(F) g_{k-i}(P/F).
\end{equation*}
\end{prop}
This generalizes the identity for simplicial polytopes
\[\sum_v g_k(P/v) = (d-k+1)g_k(P) + (k+1)g_{k+1}(P)\]
(summing over all vertices of $P$),
which is used in the proof of the upper bound theorem for
simplicial polytopes \cite{McM70}. 
The special case $d = 2k + 1$ was previously
obtained by Stenson \cite{Stenson}.

We will use Kalai's
convolution notation: if $\phi_1$ and $\phi_2$ are linear expressions
in the flag numbers of $d_1$ and $d_2$-polytopes, respectively, 
this gives an invariant $\phi\ast \psi$ of $(d_1 + d_2 + 1)$-polytopes
by
\[\phi\ast \psi(P) = \mathop{\sum_{F \le P}}_{\dim F = d_1} \phi(F)\psi(P/F);\]
it is again a linear combination of flag numbers.
 
Define $\gt_k(P) = h_k(P) - h_{k-1}(P)$ for all $k$, not just $k \le d/2$,
so by the Dehn-Sommerville relations we have
\[\gt_k(P) = \begin{cases} g_k(P) & \text{$k \le d/2$} \\
-g_{d-k+1}(P) & \text{$k \ge d/2 + 1$} \\
0 & \text{$k = (d+1)/2$\, ($d$ odd)}.
               \end{cases}\]
Let $\gt_k^d$ denote the invariant $\gt_k$ applied to $d$-polytopes.
It will be more convenient for induction to prove the following 
generalization of 
Proposition \ref{Kalai's prop}.
\begin{equation} \label{kalai's identity}
(k+1)\gt_{k+1}(P) + (d-k+1)\gt_{k}(P) = 
       \sum_{i=0}^k (i+1)\gt_i^{2i} \ast \gt^{d-2i-1}_{k-i}(P).
\end{equation}

Since this formula is linear in the flag vectors, it is enough to check
it on a basis of polytopes, i.e.\ a
collection of polytopes whose flag $f$-vectors are a 
basis for the linear span of all flag $f$-vectors.  Using the basis
given in \cite{BayBil85}, this amounts to checking \eqref{kalai's identity} 
when $d=0$, and inductively showing that if it holds for all $d-1$ polytopes $Q$,
then it holds for the cone $CQ$ and the bipyramid $BQ$. 

The case when $d$ is even and $k = d/2$ is immediate, since  
the left-hand side is $(k+1)(-g_k(P) + g_k(P)) = 0$, while each
term on the right contains $\gt_{k-i}^{2(k-i) - 1} \equiv 0$.  In particular,
the base case $d = 0$ is established. 

Next, suppose that $d > 0$ and $P = CQ$ for a $d-1$ polytope $Q$. 
By the symmetry of the $\gt_k$'s it is enough to establish 
\eqref{kalai's identity} for $0\le k < d/2$.  We have $g_k(P) = g_k(Q)$ for 
all $k$ (see \cite{Kal:new_basis}), which implies that 
$\gt_k(P) = \gt_k(Q)$ if $k \le d/2$.  Thus
if $k \le d/2 - 1$ the left side of \eqref{kalai's identity} becomes
\[(k+1)\gt_{k+1}(Q) + (d-k+1)\gt_k(Q) = \gt_k(Q) +
\sum_{i=0}^k (i+1)g_i^{2i} \ast \gt_{k-i}^{d-2i-2}(Q),
\]
by the inductive hypothesis.  On the other hand, the nonempty faces
of $CQ$ are either nonempty faces of $Q$ or cones over faces
(possibly empty).  Since $g_k(CF) = 0$ if $\dim(CF)=d$, the
only $2i$-faces of $P$ for which $g_i \ne 0$ are the apex 
$C\emptyset$ and faces of $Q$.  Thus the right side of \eqref{kalai's identity} becomes
\[\gt_k(Q) +
\sum_{i=0}^k (i+1)g_i^{2i} \ast \gt_{k-i}^{d-2i-2}(Q),
\]
as required.

If $d$ is odd and $k = (d-1)/2$, we have $\gt_k(P) = \gt_k(Q)$, but
$\gt_{k+1}(P) = 0$, so the left side of \eqref{kalai's identity} is
just $(k+2)\gt_k(Q)$.  Since for every face $F$ of $P$ other than 
the apex $C\emptyset$ and the base $Q$ either $F$ or $P/F$ is a 
cone over a nonempty polytope, the right side is 
\[\sum_{i=0}^k (i+1)g_i^{2i} \ast \gt_{k-i}^{2(k-i)}(P) = 
g_0(C\emptyset)g_k(Q) + (k+1)g_k(Q)g_0(P/Q) = (k+2)g_k(Q).\]
Therefore \eqref{kalai's identity} holds for $P = CQ$.

Now suppose that $P = BQ$, the bipyramid over a $(d-1)$-polytope $Q$.
Then $h_{P}(t) = (t+1)h_Q(t)$ (see \cite{Kal:new_basis}), so we have 
\[\gt_k(P) = \gt_k(Q) + \gt_{k-1}(Q)\]
for all $0 \le k \le d+1$.  Thus the
left side of our identity is
\[ (k+1)(\gt_{k+1}(Q) + \gt_{k}(Q)) + (d - k + 1)(\gt_k(Q) + \gt_{k-1}(Q)) 
\]
\begin{equation} \label{eq1}
 = 2\gt_k(Q) + \sum_{i=0}^k (i+1)g_i^{2i} \ast \gt_{k-i}^{d-2i-2}(Q) + 
\sum_{i=0}^{k-1} (i+1)g_i^{2i} \ast \gt_{k-1-i}^{d-2i-2}(Q),
\end{equation}
by the inductive hypothesis.

On the other hand, faces of $P = BQ$ are either faces of $Q$ other than 
$Q$ itself or cones over faces of $Q$ (possibly empty).  The 
only $2i$-faces $F \le BQ$ for which 
$g_i^{2i}(F)$ can be nonzero are the two apexes, for which
$g_0(F)\gt_k(BQ/F) = \gt_k(Q)$, and faces $F \le Q$, $F \ne Q$,
for which 
\[g_i(F)g_{k-i}((BQ)/F) = g_i(F)\gt_{k-i}(B(Q/F)) = 
g_i(F)[\gt_{k-i}(Q/F) + \gt_{k-i-1}(Q/F)].\]
Substituting these into the right side of \eqref{kalai's identity}  
gives \eqref{eq1}. 

\bibliographystyle{abbrv}
\bibliography{expos}
\end{document}